\documentclass[11pt]{article}
\usepackage{amsfonts}

\usepackage{amsmath}
\usepackage{amssymb}
\usepackage{latexsym}
\newtheorem{Theorem}{Theorem}[section]
\newtheorem{Lemma}{Lemma}[section]
\newtheorem{Remark}{Remark}[section]
\newtheorem{cor}{Corollary}[section]

\newcommand\qed{\hfill$\Box$}
\newcommand{\beq}{\begin{eqnarray}}
\newcommand{\eeq}{\end{eqnarray}}
\newcommand{\beqno}{\begin{eqnarray*}}
\newcommand{\eeqno}{\end{eqnarray*}}
\newcommand{\be}{\begin{equation}}
\newcommand{\ee}{\end{equation}}
\newcommand{\ideq}{{\lower .5ex
\hbox{$\>\>\stackrel{\triangle}{=}\>\>$} }}

\def \d2{\Delta_{+}^2}
\begin{document}

\newcommand{\D}{\displaystyle}
\title{\bf Global Classical Large Solutions to Navier-Stokes Equations for Viscous
Compressible and Heat Conducting Fluids with Vacuum}
\author{
\begin{tabular}{cc}
& H{\sc uanyao} W{\sc en}$^{1,2}$,  \ \ \ C{\sc hangjiang} Z{\sc
hu$^{2}$}\thanks{Corresponding author.\ \ Email:
cjzhu@mail.ccnu.edu.cn} \\[4mm]
$^1$ &  School of Mathematical Sciences\\
& South China Normal University,  Guangzhou  510631,  P.R.
China\\[2mm]
$^2$ &  The Hubei Key Laboratory of Mathematical Physics \\
& School of Mathematics and Statistics
\\
& Huazhong Normal University, Wuhan 430079, P.R. China\\
\end{tabular}
}
\date{}

\maketitle

\begin{abstract} In this paper, we consider the 1D Navier-Stokes equations for viscous
compressible and heat conducting fluids (i.e., the full
Navier-Stokes equations). We get a unique global classical solution
to the equations with large initial data and vacuum. Because of the
strong nonlinearity and degeneration of the equations brought by the
temperature equation and by vanishing of density (i.e., appearance
of vacuum) respectively, to our best knowledge, there are only two
results until now about global existence of solutions to the full
Navier-Stokes equations with special pressure, viscosity and heat
conductivity when vacuum appears (see \cite{Feireisl-book} where the
viscosity $ \mu=$const and the so-called {\em variational} solutions
were obtained, and see \cite{Bresch-Desjardins} where the viscosity
$ \mu=\mu(\rho)$ degenerated when the density vanishes and the
global weak solutions were got). It is open whether the global
strong or classical solutions exist. By applying our ideas which
were used in our former paper \cite{Ding-Wen-Zhu} to get
$H^3-$estimates of $u$ and $\theta$ (see Lemma \ref{non-le:3.14},
Lemma \ref{non-le:3.15}, Lemma \ref{non-rle:3.12} and the
corresponding corollaries), we get the existence and uniqueness of
the global classical solutions (see Theorem \ref{non-rth:1.1}). In
fact, the existence of strong solutions would be done obviously by
our estimates if the regularity of the initial data is assumed to be
weaker. Like \cite{Ding-Wen-Zhu}, we get $H^4-$regularity of $\rho$
and $u$ (see Theorem \ref{non-th:1.1}). We do not get further
regularity of $\theta$ such as $H^4-$regularity, because of the
degeneration and strong nonlinearity brought by vacuum and the term
$(\mu uu_x)_x$ in the temperature equation. This can be viewed the
first result on global classical solutions to the 1D Navier-Stokes
equations for viscous compressible and heat conducting fluids which
may be large initial data and contain vacuum.
\end{abstract}


%\bigbreak
\section*{Contents}

1. Introduction \dotfill 2

\noindent 2. Preliminaries \dotfill 7

\noindent 3.  Proof of Theorem \ref{non-rth:1.1} \dotfill 9

\noindent 4. Proof of Theorem \ref{non-th:1.1} \dotfill 26

\noindent References \dotfill 36

\vspace{4mm}

\setcounter{section}{0} \setcounter{equation}{0}
\section{Introduction}
In this paper, we consider the Navier-Stokes equations for viscous
compressible and heat conducting fluids (i.e. the full Navier-Stokes equations). The model, describing for
instance the motion of gas, plays an important role in applied
physics. Mathematically, the model in one dimension can be written
as follows in sense of Eulerian coordinates: \beq \label{non-1.2}
\begin{cases}
\rho_t+(\rho u)_x=0,\ \rho\ge0,\\
(\rho u)_t+(\rho u^2)_x+P_x=(\mu u_x)_x,\\
(\rho E)_t+(\rho u E)_x+(Pu)_x=(\mu uu_x)_x+(\kappa \theta_x)_x,
\end{cases} \eeq for $(x,t)\in (0,1)\times (0,+\infty)$.
Here $\rho=\rho(x,t)$, $u=u(x, t)$, $P=P(\rho,\theta)$, $E$,
$\theta$ and $\kappa=\kappa(\rho, \theta)$ denote the density,
velocity, pressure, total energy, absolute temperature and
coefficient of heat conduction, respectively. The total energy
$E=e+\frac{1}{2}u^2$, where $e$ is the internal energy.
 $\mu>0$ is the coefficient of viscosity.
 $P$ and $e$ satisfy the second principle
of thermodynamics: \be\label{non-1.1} P=\rho^2\frac{\partial
e}{\partial\rho}+\theta\frac{\partial P}{\partial \theta}. \ee In
the present paper, we consider the initial and boundary conditions:
\be\label{non-1.3}(\rho, \ u,\ \theta)\big|_{t=0}=(\rho_0,\ u_0,\
\theta_0)(x) \ \ {\rm{in}} \ \ [0,1], \ee and\be\label{non-1.4}
 (u,\ \theta_x)\big|_{x=0,1}=0, \ t\geq0.
\ee
Since the model is important, lots of works on the
existence,
 uniqueness, regularity and asymptotic behavior of the solutions
 were done during the last five decades. While, because of the stronger nonlinearity in
(\ref{non-1.2}) compared with the Navier-Stokes equations for
isentropic fluids (no temperature equation), many known mathematical
results mainly focused on the absence of vacuum (vacuum means
$\rho=0$), refer for instance to \cite{Itaya, Jiang1, 6,
Kazhikhov-Shelukhi, 2, 3, 1} for classical solutions. More
precisely, the local classical solutions to the Navier-Stokes
equations with heat-conducting fluid in H\"older spaces was obtained
respectively by Itaya in \cite{Itaya} for Cauchy problem and by Tani
in
 \cite{1} for IBVP with $\inf\rho_0>0$, where the spatial dimension $N=3$. Using
delicate energy methods in Sobolev spaces, Matsumura and Nishida
showed in \cite{2,3} that the global classical solutions exist
provided that the initial data is small in some sense and away from
vacuum and the spatial dimension $N=3$. For large initial data and
dimension $N=1$, Kazhikhov, Shelukhi in \cite{Kazhikhov-Shelukhi}
(for polytropic perfect gas with $\mu,\kappa=$const.) and Kawohl in
\cite{6} (for real gas with $\kappa=\kappa(\rho,\theta),\
\mu=\mathrm{const.}$) respectively got global classical solutions to
(\ref{non-1.2}) in Lagrangian coordinates with boundary condition
(\ref{non-1.4}) and $\inf\rho_0>0$. The internal energy $e$ and the
coefficient of heat conduction $\kappa$ in \cite{6} satisfy the
following assumptions for $\rho\le\overline{\varrho}$ and
$\theta\ge0$ (we translate these conditions in Eulerian coordinates)
\be\label{non-r1.5}\begin{cases} e(\rho,0)\ge0,\ \
\nu(1+\theta^r)\le\partial_\theta e(\rho,\theta)\le
N(\overline{\varrho})(1+\theta^r),\\
\kappa_0(1+\theta^q)\le\kappa(\rho,\theta)\le\kappa_1(1+\theta^q),\\
|\partial_\rho\kappa(\rho, \theta)|+|\partial_{\rho\rho}\kappa(\rho,
\theta)|\le\kappa_1(1+\theta^q),
\end{cases}\ee
where $r\in[0,1]$, $q\ge2+2r$, and $\nu$, $N(\overline{\varrho})$,
$\kappa_0$ and $\kappa_1$ are positive constants. For the
perfect gas in the domain exterior to a ball in $\mathbb{R}^N$ ($N =
2$ or $3$) with $\mu,\kappa=$const., Jiang in \cite{Jiang1} got the
existence of global classical spherically symmetric large solutions
in H\"older spaces.

In fact, Kawohl in \cite{6} also considered the case of
$\mu=\mu(\rho)$ for another boundary condition with $\inf\rho_0>0$,
where $0<\underline{\mu}_0\le\mu(\rho)\le\overline{\mu}_0$ for any
$\rho\ge0$ and $\underline{\mu}_0$ and $\overline{\mu}_0$ are
positive constant. This result was generalized to the case
$\mu(\rho)=\rho^\alpha$ by Jiang in \cite{18} for
$\alpha\in(0,\frac{1}{4})$, and by Qin, Yao in \cite{19} for
$\alpha\in(0,\frac{1}{2})$, respectively.

On the existence, asymptotic behavior of the weak solutions for full
Navier-Stokes equations (including the temperature equation) with
$\inf\rho_0>0$, please refer for instance to \cite{ Jiang2, Jiang3,
Jiang-Zhang:weak solutions} for weak solutions in 1D and for
spherically symmetric weak solutions in bounded annular domains in
$\mathbb{R}^N$ ($N=2$, $3$), and refer to \cite{Feireisl1} for {\em
variational} solutions in a bounded domain in $\mathbb{R}^N$ ($N=2$,
$3$).

In the presence of vacuum (i.e. $\rho$ may vanish), to our best
knowledge, the mathematical results about global well-posedness of
the full Navier-Stokes equations are usually limited to the
existence of weak solutions with special pressure, viscosity and
heat conductivity (see \cite{Bresch-Desjardins, Feireisl-book}).
More precisely, Feireisl in \cite{Feireisl-book} got the existence
of so-called {\em variational} solutions in dimension $N\ge 2$. The
temperature equation in \cite{Feireisl-book} is satisfied only as an
inequality. Anyway, this work in \cite{Feireisl-book} is the very
first attempt towards the existence of weak solutions to the full
compressible Navier-Stokes equations in higher dimensions, where the
viscosity $\mu$ is {\em constant} and
\be\label{non-1.6}\begin{cases} \kappa=\kappa(\theta)\in
C^2[0,\infty),\
\underline{\kappa}(1+\theta^a)\le\kappa(\theta)\le\overline{\kappa}(1+\theta^a)\
\ \mathrm{for}\ \mathrm{all}\ \theta\ge0,
\\ P=P(\rho,\theta)=\mathcal {P}_e(\rho)+\theta\mathcal{P}_\theta(\rho)\ \
\mathrm{for}\ \mathrm{all}\ \rho\ge0\ \mathrm{and}\ \theta\ge0,\\ \mathcal{P}_e,
\mathcal{P}_\theta\in C[0,\infty)\cap C^1(0,\infty);\ \mathcal{P}_e(0)=0,\ \mathcal{P}_\theta(0)=0,
\\ \mathcal{P}_e^\prime(\rho)\ge a_1\rho^{\overline{\gamma}-1}-b_1\ \ \mathrm{for}\
\mathrm{all}\ \rho>0;\ \mathcal{P}_e(\rho)\le a_2\rho^{\overline{\gamma}}+b_1\ \
\mathrm{for}\ \mathrm{all}\ \rho\ge0,\\
\mathcal{P}_\theta \ \ \mathrm{is} \ \ \mbox{non-decreasing}\
\mathrm{in}\ [0,\infty);\ \mathcal{P}_\theta(\rho)\le
a_3(1+\rho^\Gamma)\ \ \mathrm{for}\ \mathrm{all}\ \rho\ge0,
\end{cases}\ee
where
$\Gamma<\frac{\overline{\gamma}}{2}$ if $N=2$ and $\Gamma=\frac{\overline{\gamma}}{N}$
 for $N\ge3$; $a\ge2$, $\overline{\gamma}>1$, and $a_1$, $a_2$, $a_3$, $b_1$, $\underline{\kappa}$
 and $\overline{\kappa}$ are positive constants. Note that the perfect gas
equation of state (i.e. $P=R\rho\theta$ for some constant $R>0$) is not involved in (\ref{non-1.6}).
In order that the equations are satisfied as equalities in the sense of distribution, Bresch and Desjardins in
\cite{Bresch-Desjardins} proposed some different assumptions from \cite{Feireisl-book}, and obtained the existence of global weak
solutions to the full Navier-Stokes equations with large initial data in $\mathbb{T}^3$
or $\mathbb{R}^3$. In \cite{Bresch-Desjardins}, the viscosity $\mu=\mu(\rho)$ may vanish when vacuum
appears, and $\kappa$, $P$ and $e$ are assumed to
satisfy\be\label{non-1.7}\begin{cases}
\kappa(\rho,\theta)=\kappa_0(\rho,\theta)(\rho+1)(\theta^a+1),\\
P=r\rho\theta+p_c(\rho), \\
e=C_\upsilon\theta+e_c(\rho),
\end{cases}\ee
where $a\ge2$, $r$ and $C_\upsilon$ are two positive constants,
$p_c(\rho)=O(\rho^{-\ell})$ and $e_c(\rho)=O(\rho^{-\ell-1})$ (for
some $\ell>1$) when $\rho$ is small enough, and
$\kappa_0(\rho,\theta)$ is assumed to satisfy
$$\underline{c}_0\le\kappa_0(\rho,\theta)\le\frac{1}{\underline{c}_0},$$
for $\underline{c}_0>0$. We have to mention that the smooth solutions in
$C^1\left([0,\infty);H^d(\mathbb{R}^N)\right)$ ($d>2+[\frac{N}{2}]$) would blow up when the initial density is of nontrivial
compact support (see \cite{25}). On the local existence and uniqueness of strong solutions in $\mathbb{R}^3$, please refer to
\cite{cho-Kim: perfect gas} for the perfect gas with
$\mu,\kappa=$const.

It is still open whether global strong (or classical) solutions exist when
vacuum appears (i.e., the density may vanish). Our main concern here
is to show the existence and uniqueness of global classical
solutions to (\ref{non-1.2})-(\ref{non-1.4}) with vacuum and large initial data. In fact,
the existence of the strong solutions to this problem is obvious if
the regularity of initial data is assumed to be weaker.

For compressible isentropic Navier-Stokes equations (i.e. no
temperature equation), there are so many results about the
well-posedness and asymptotic behaviors of the solutions when vacuum
appears. Refer to \cite{Feireisl2, Jiang-Zhang1, Lions2, Luo} and
\cite{Guo-Zhu, 13, 23, 15, 24, Zhu} for global weak solutions with
constant viscosity and with density-dependent viscosity,
respectively. Refer to \cite{9, 10} and \cite{7, 11, 8, salvi} for
global strong solutions and for local strong (classical) solutions
with constant viscosity, respectively. Recently, Huang, Li, Xin in
\cite{Huang-Li-Xin} and Ding, Wen, Yao, Zhu in \cite{Ding-Wen-Zhu,
Ding-Wen-Yao-Zhu} independently got existence and uniqueness of
global classical solutions, where the initial energy in
\cite{Huang-Li-Xin} is assumed to be small in $\mathbb{R}^3$ and
$\rho-\widetilde{\rho}\in C\left([0,T]; H^3(\mathbb{R}^3)\right)$,
$u\in C\left([0,T];D^1(\mathbb{R}^3)\cap
D^3(\mathbb{R}^3)\right)\cap L^\infty\left([\tau, T];
D^4(\mathbb{R}^3)\right)$ (for $\tau>0$) which generalized the
results in \cite{11}, and the initial data in \cite{Ding-Wen-Zhu,
Ding-Wen-Yao-Zhu} could be large for dimension $N=1$ and could be
large but spherically symmetric for $N\ge2$, and $(\rho,u)\in
C([0,T]; H^4(I))$ ($I$ is bounded in \cite{Ding-Wen-Zhu}, and is
bounded or an exterior domain in
\cite{Ding-Wen-Yao-Zhu}).\\[1.5mm]

We would like to give some notations which will be used throughout
the paper.

 \noindent{\bf Notations:}

(1) $I=[0,1]$, $\partial I=\{0,1\}$, $Q_T=I\times[0,T]$ for $T>0$.

(2) For $p\in[1, \infty]$, $L^p=L^p(I)$ denotes the $L^p$ space with
the norm $\|\cdot\|_{L^p}$. For $k\ge 1$ and $p\in[1, \infty]$,
$W^{k,p}=W^{k,p}(I)$ denotes the Sobolev space, whose norm is
denoted as $\|\cdot\|_{W^{k,p}}$, $H^k=W^{k,2}(I)$.

 (3) For an integer $k\ge 0$ and $0<\alpha<1$, let $C^{k+\alpha}(I)$ denote
the Schauder space of functions on $I$, whose $k$th order derivative
is H\"older continuous with exponents $\alpha$, with the norm
$\|\cdot\|_{C^{k+\alpha}}$.\\[1.5mm]

In this paper, our assumptions are the following:

($A_1$): $\int_I\rho_0>0$.

($A_2$): $\mu=\mathrm{const.}>0$, $e=C_0Q(\theta)+e_c(\rho)$,
$P=\rho Q(\theta)+P_c(\rho)$, $\kappa=\kappa(\theta)$, for some
constant $C_0>0$.

($A_3$): $P_c(\rho)\ge0$, $e_c(\rho)\ge0$, for $\rho\ge0$; $P_c\in
C^2[0,\infty)$; $\rho| \frac{\partial e_c}{\partial\rho}|\le C_1
e_c(\rho)$, for some constant $C_1>0$.

$(A_4):$ $Q(\cdot)\in C^2[0,\infty)$ satisfies
$$\begin{cases}C_2\left(\beta+(1-\beta)\theta+\theta^{1+r}\right)\le
Q(\theta)\le
C_3\left(\beta+(1-\beta)\theta+\theta^{1+r}\right),\\
C_4(1+\theta^r)\le Q'(\theta)\le C_5(1+\theta^r) ,\end{cases}\ \ \ \
\ \  \ \ \ \ \ \ \
$$
for some constants $C_i>0$ ($i= 2$, $3$, $4$, $5$) and $r\ge 0$,
$\beta=0$
or $1$.\\

($A_5$): $\kappa\in C^2[0,\infty)$ satisfies
 $$C_6(1+\theta^q)\le \kappa(\theta)\le
C_7(1+\theta^q),$$ for
$q\ge2+2r$, and some constants $C_i>0$ ($i=6$, $7$).\\

($A_6$): $Q, P_c\in C^4[0,\infty)$, and $\kappa$ satisfies
$$
|\partial_\theta^3\kappa(\theta)|\le C_8(1+\theta^{q-3}),$$ for
$\theta>0$ and some constant $C_8>0$.
\begin{Remark}\label{non-remark1.1}
$(i)$ $(A_1)$ is needed to get the upper bounds of $\theta$ and
$\theta_t$ in terms of some norms by using mass conservation, Lemma
\ref{non-le:2.1} and Corollary
\ref{non-cor:2.1}.\\[1.5mm]

$(ii)$ The case for the perfect gas (i.e. $P=R\rho\theta$,
$e=C_\nu\theta$ for constants $R>0$ and $C_\nu>0$) is
involved in the above assumptions. \\[1.5mm]

$(iii)$ As it mentioned in \cite{6}, the restriction on $\mu$
($\mu$=const., see other restrictions on $\kappa$ and e in
(\ref{non-r1.5})) is not physically motivated. Physically, it seems
more importantly that the state functions e, $\mu$ and $\kappa$
usually depend on both $\rho$ and $\theta$. Particularly, the
internal energy e grows as $\theta^{1+r}$ with $r\approx 0.5$, the
conductivity $\kappa$ grows as $\theta^q$ with $4.5\le q\le5.5$ and
viscosity $\mu$ increases like $\theta^p$ with $0.5\le p\le0.8$ (see
\cite{6, 19} and references therein). Because of mathematical
technique, in the present paper, we assume $\mu=$const. and
$\kappa=\kappa(\theta)$ as in \cite{Feireisl-book} (see
(\ref{non-1.6})). From ($A_2$)--($A_5$), we know that e and $\kappa$
grow respectively as $\theta^{1+r}$ and $\theta^q$, where $q$ can be
taken as $q\in[4.5,5.5]$, and r can be taken as $r=0.5$ if we
consider $\theta>0$.\\[1.5mm]

$(iv)$  The restriction on $q$ in ($A_5$) (i.e. $q\ge 2+2r$) is same
as (\ref{non-r1.5}), and is the same as (\ref{non-1.6}) and
(\ref{non-1.7}) when we take $r=0$. This assumption plays an
important role in the analysis.
\end{Remark}

\noindent{\bf Main results:}

\begin{Theorem}\label{non-rth:1.1}
In addition to $(A_1)$-$(A_5)$, we assume $\rho_0\geq0$,
 $\rho_0\in H^2$, $(\sqrt{\rho_0})_x\in L^\infty$, $u_0\in H^3\cap H_0^1$,
$\theta_0\in H^3$, $\partial_x\theta_0|_{x=0,1}=0$, and that the
following compatible conditions are valid: \beq\label{non-1.8}
\begin{cases}\mu u_{0xx}-[P(\rho_0,\
\theta_0)]_x=\sqrt{\rho}_0g_1,\\
\left[\kappa(\theta_0)
\theta_{0x}\right]_x+\mu|u_{0x}|^2=\sqrt{\rho_0}\ g_2,\ x\in I,
\end{cases}\eeq for some $g_1,g_2\in L^2$, and $\left(\sqrt{\rho_0}g_1\right)_x, \left(\sqrt{\rho_0}g_2\right)_x\in L^2$. Then for any $T>0$ there exists
a unique global solution $(\rho,\ u,\ \theta)$ to
(\ref{non-1.2})-(\ref{non-1.4}) such that \beqno &\rho\in
C([0,T];H^2),\ \rho_t\in C([0,T];H^1),\ \sqrt{\rho}\in
W^{1,\infty}(Q_T), &\\& u\in L^\infty([0,T];H^3),\ \sqrt{\rho}u_t\in
L^\infty([0,T]; L^2),&\\&\rho u_t\in L^\infty([0,T]; H^1),\ \ u_t\in
L^2([0,T]; H_0^1),\ \ \sqrt{\rho}e_t\in L^\infty([0,T];
L^2),&\\&\rho e_t\in L^\infty([0,T]; H^1), \ \ \theta\in
L^\infty([0,T];
 H^3),\ \ \theta_t\in L^2([0,T]; H^1).&
\eeqno
\end{Theorem}
\begin{Remark}
(i)\  (\ref{non-1.8}) was proposed by Cho and Kim in \cite{cho-Kim:
perfect gas} to get $\mathrm{local}$ $H^2$-{\rm regularity} of $u$
and $\theta$ for $\mathrm{the\ polytropic\ perfect\ gas}$. The
detailed reasons why such conditions were needed can be found in
\cite{cho-Kim: perfect gas}. Roughly speaking, $g_1$ and $g_2$ are
equivalent to $\sqrt{\rho}u_t$ and $\sqrt{\rho}e_t$  at $t=0$,
respectively.

(ii)\ \ \ By the Sobolev embedding theorems (cf. \cite{26}) and
Lemma \ref{non-le:2.3}, we know from Theorem \ref{non-rth:1.1}
\beqno &&\rho\in C\left([0,T]; C^{1+\frac{1}{2}}(I)\right)\cap
C^1\left([0,T];C^{\frac{1}{2}}(I)\right),\\&& u\in  C\left([0,T];
C^{2+\sigma}(I)\right),\ (\rho u)_t\in  C\left([0,T];
C^{\sigma}(I)\right),\\&& \theta\in  C\left([0,T];
C^{2+\sigma}(I)\right),\ (\rho e)_t\in  C\left([0,T];
C^{\sigma}(I)\right), \eeqno for any $T>0$ and
$\sigma\in(0,\frac{1}{2})$. This implies $(\rho, u, \theta)$ is the
classical solution to (\ref{non-1.2})-(\ref{non-1.4}).
\end{Remark}
\begin{Theorem}\label{non-th:1.1}
In addition to $(A_1)$-$(A_6)$, we assume $\rho_0\geq0$,
 $\rho_0\in H^4$, $(\sqrt{\rho_0})_x\in L^\infty$, $u_0\in H^4\cap H_0^1$,
$\theta_0\in H^3$, $\partial_x\theta_0|_{x=0,1}=0$, $q>2+2r$, and
that the following compatible conditions are valid:
\beq\label{non-1.5}
\begin{cases}\mu u_{0xx}-[P(\rho_0,\
\theta_0)]_x=\rho_0g_3,\\
\left[\kappa(\theta_0)
\theta_{0x}\right]_x+\mu|u_{0x}|^2=\sqrt{\rho_0}\ g_2,\ x\in I,
\end{cases}\eeq for some $g_3\in H_0^1$, $\left(\sqrt{\rho_0}\partial_xg_3\right)_x\in L^2$, and $g_2,\left(\sqrt{\rho_0}g_2\right)_x\in L^2$. Then for any $T>0$ there exists
a unique global solution $(\rho,\ u,\ \theta)$ to
(\ref{non-1.2})-(\ref{non-1.4}) satisfying: \beqno &\rho\in
C([0,T];H^4),\ \rho_t\in C([0,T];H^3),\ \sqrt{\rho}\in
W^{1,\infty}(Q_T), &\\& u\in C([0,T];H^4)\cap L^2([0,T];H^5),\
u_t\in L^\infty([0,T];H_0^1)\cap L^2([0,T];H^3),&\\&(\rho u)_t\in
C([0,T]; H^2),\ \ \sqrt{\rho}u_{xxt}\in L^\infty([0,T];L^2),\ \
\sqrt{\rho}e_t\in L^\infty([0,T]; L^2),&\\&(\rho e)_t\in
L^\infty([0,T]; H^1), \ \ \theta\in L^\infty([0,T];
 H^3)\cap L^2([0,T]; H^4),\ \ \theta_t\in L^2([0,T]; H^1).&
\eeqno
\end{Theorem}
\begin{Remark}
(i)\ (\ref{non-1.5})$_1$ was proposed by Cho and Kim in \cite{11}
where they consider the local existence of classical solutions for
$\mathrm{isentropic\ fluids}$ (no temperature equation). Roughly
speaking, $g_3$ is equivalent to $u_t$ at $t=0$.

(ii)\ We could not get\ \  $\theta\in C([0,T]; H^4)$ (or $L^\infty([0,T]; H^4)$) even if
(\ref{non-1.5})$_2$ is changed similarly to (\ref{non-1.5})$_1$,
because of the strong nonlinearity and degeneration brought by $(\mu
uu_x)_x$ in the temperature equation and the appearance of vacuum,
respectively.

(iii)\ Using ideas of Cho and Kim in \cite{11}, we can also get
$$
u\in L^\infty\left([\tau, T];H^4\right),\ \theta\in
L^\infty\left([\tau, T];H^3\right),
$$ for $\tau>0$. If we can obtain our estimates in higher
dimensions, it will be useful to investigate the local (global)
existence of classical solutions to the full Navier-Stokes equations
(including the temperature equation) in $\mathbb{R}^N$ ($N\ge2$).
For example, to guarantee (\ref{non-1.4}) in $\mathbb{R}^N$ ($N=2$
or $3$) is valid for all $t\ge0$, it is necessary to get\ \
$\theta\in L^\infty\left([0, T];H^3\right)$. We will consider these
problems in the near future.
\end{Remark}
The constants $C_0$ in $(A_2)$ and the viscosity $\mu$ don't play
any role in the analysis, we assume henceforth that $C_0=1$ and
$\mu=1$ for simplicity.

The rest of the paper is organized as follows. In Section 2, we
present some useful lemmas which will be used in the next sections.
In Section 3, we prove Theorem \ref{non-rth:1.1}
by giving the initial density and the initial temperature a lower
bound $\delta>0$, getting a sequence of approximate solutions to
(\ref{non-1.2})-(\ref{non-1.4}), and taking $\delta\rightarrow0^+$ after
making some estimates uniformly for $\delta$. More precisely, based on Lemma \ref{non-le:2.1} and the one-dimensional properties of the equations,
we get $H^2-$estimates of the solutions. Using our ideas in \cite{Ding-Wen-Zhu, Ding-Wen-Yao-Zhu}, we obtain $H^3-$estimates of $u$ and $\theta$. In Section 4, using
the similar arguments as in Section 3, we prove Theorem
\ref{non-th:1.1}.

\setcounter{section}{1} \setcounter{equation}{0}
\section{Preliminaries}
\begin{Lemma}\label{non-le:2.1}
Let $\Omega=[\overline{a}, \overline{b}]$ be a bounded domain in $\mathbb{R}$, and
$\rho$ be a non-negative function such that
$$
0<M\le\int_\Omega\rho\le K,
$$ for constants $M>0$ and $K>0$. Then
$$ \|v\|_{L^\infty(\Omega)}\le \frac{K}{M}\|v_x\|_{L^1(\Omega)}+\frac{1}{M}\left|\int_\Omega\rho v\right|,
$$
 for any $v\in H^1(\Omega)$.
\end{Lemma}
\noindent{\it Proof.} For any $x\in\Omega$, we have \beqno
|v(x)|&\le&\frac{1}{M}\left|v(x)\int_\Omega\rho(y)dy\right| \\
[4mm]
&\le &\frac{1}{M}\left|\int_\Omega
v(x)\rho(y)dy-\int_\Omega\rho(y)v(y)dy\right|+\frac{1}{M}\left|\int_\Omega\rho(y)v(y)dy\right|
\\[4mm]
& \le & \frac{1}{M}\left|\int_\Omega
\int_y^xv_\xi(\xi)d\xi\rho(y)dy\right|+\frac{1}{M}\left|\int_\Omega\rho(y)v(y)dy\right|
\\[4mm]
& \le &
\frac{K}{M}\|v_x\|_{L^1(\Omega)}+\frac{1}{M}\left|\int_\Omega\rho(y)v(y)dy\right|.
 \eeqno
\qed
\begin{Remark}
The version of higher dimensions for Lemma \ref{non-le:2.1} can be
found in \cite{Feireisl1} or \cite{Feireisl-book}.
\end{Remark}

\begin{cor}\label{non-cor:2.1}
Consider the same conditions in Lemma \ref{non-le:2.1}, and in
addition assume $\Omega=I$, and
$$
\|\rho v\|_{L^1(I)}\le \overline{c}.
$$
Then for any $l>0$, there exists a positive constant $C=C(M, K, l,
\overline{c})$ such that $$\|v^l\|_{L^\infty(I)}\le
C\|(v^l)_x\|_{L^2(I)}+C,$$ for any $v^l\in H^1(I)$.
\end{cor}
\noindent{\it Proof.} By Lemma \ref{non-le:2.1}, we have \beqno
\|v^l\|_{L^\infty(I)}\le C\|(v^l)_x\|_{L^2(I)}+C\int_I\rho |v^l|.
\eeqno Case 1: $l\in (0, 1]$.

In this case, we use the Young inequality to get
 \beqno \|v^l\|_{L^\infty(I)}&\le&C\|(v^l)_x\|_{L^2(I)}+C\int_I\rho
|v|+C\int_I\rho+C\\&\le&C\|(v^l)_x\|_{L^2(I)}+C.\eeqno
 Case
2: $l\in(1, \infty)$.

In the case, we use the Young inequality again to get
 \beqno \|v^l\|_{L^\infty(I)}&\le&C\|(v^l)_x\|_{L^2(I)}+C\|v^{l-1}\|_{L^\infty(I)}\int_I\rho
|v|\\&\le&C\|(v^l)_x\|_{L^2(I)}+\frac{1}{2}\|v^l\|_{L^\infty(I)}+C.\eeqno
This gives
 \beqno \|v^l\|_{L^\infty(I)}\le C\|(v^l)_x\|_{L^2(I)}+C.\eeqno
\qed
\begin{Lemma}\label{non-le:2.2}
For any $v\in H^1_0(I)$, we have
$$
\|v\|_{L^\infty(I)}\le \|v_x\|_{L^1}.
$$
\end{Lemma}
\noindent{\it Proof.} Since $v(0)=0$, we have for any $x\in I$
\beqno
|v(x)|=|v(x)-v(0)|=\left|\int_0^xv_x\right|\le\|v_x\|_{L^1(I)}.
\eeqno This completes the proof. \qed

\begin{Lemma}\label{non-le:2.3} (\cite{Simon}).
Assume $X\subset E\subset Y$ are Banach spaces and
$X\hookrightarrow\hookrightarrow E$. Then the following imbedding
are compact: $$\D(i)\ \  \left\{\varphi:\varphi\in L^q(0,T; X),
\frac{\partial\varphi}{\partial t}\in
L^1(0,T;Y)\right\}\hookrightarrow\hookrightarrow L^q(0,T; E),\ \
{\rm if}\ \  1\leq q\leq\infty; $$ $$(ii)\ \
\left\{\varphi:\varphi\in L^\infty(0,T; X),
\frac{\partial\varphi}{\partial t}\in
L^r(0,T;Y)\right\}\hookrightarrow\hookrightarrow C([0,T]; E),\ \
{\rm if}\ \   1< r\leq\infty. $$
\end{Lemma}
 \setcounter{section}{2}
\setcounter{equation}{0}
\section{ \ Proof of Theorem \ref{non-rth:1.1}}
In this section, we get a global solution to
(\ref{non-1.2})-(\ref{non-1.4}) with initial density and initial
temperature having a respectively lower bound $\delta>0$ by using
some {\it a priori} estimates of the solutions based on the local
existence. Theorem \ref{non-rth:1.1} would be got after we make some
{\it a priori} estimates uniformly for $\delta$ and take
$\delta\rightarrow0^+$.

Denote $\rho_0^\delta=\rho_0+\delta$ and
$\theta_0^\delta=\theta_0+\delta$ for $\delta\in(0,1)$. Throughout
this section, we denote $c$ to be a generic constant depending on
$\rho_0$, $u_0$, $\theta_0$, $T$ and some other known constants but
independent of $\delta$ for any $\delta\in(0,1)$.\\[0.4mm]

 Before
proving Theorem \ref{non-rth:1.1}, we need the following auxiliary
theorem.
\begin{Theorem}\label{non-th:3.1} Consider the same assumptions as in Theorem \ref{non-rth:1.1}. Then for any $T>0$ and $\delta\in(0,1)$ there exists a unique global solution $(\rho,
u, \theta)$ to (\ref{non-1.2})-(\ref{non-1.4}) with initial data
replaced by ($\rho_0^\delta,u_0,\theta_0^\delta$), such that \beqno
&\rho\in C([0,T];H^2), \ \ \ \rho_t\in C([0,T];H^1), \ \ \
\rho_{tt}\in L^2([0,T];L^2), \ \rho\ge\frac{\D\delta}{c}>0, &\\&
u\in C([0,T];H^3\cap H^1_0),\ u_t\in C([0,T];H^1)\cap
L^2([0,T];H^2), \ \ \ u_{tt}\in
 L^2([0,T];L^2),\ &\\& \theta\ge
c_\delta>0,\ \theta\in C([0,T]; H^3),\ \ \theta_t\in C([0,T];
H^1)\cap L^2([0,T]; H^2),\ \ \theta_{tt}\in L^2([0,T]; L^2), &
\eeqno where $c_\delta$ is a constant depending on $\delta$, but
independent of $u$.
\end{Theorem}

{\noindent\bf Proof of Theorem \ref{non-th:3.1}:}

The local solutions as in Theorem \ref{non-th:3.1} can be obtained
by the successive approximations like in \cite{cho-Kim: perfect
gas}. We omit it here for brevity. The regularities guarantee the
uniqueness (refer for instance to \cite{cho-Kim: perfect gas}).
Based on it, Theorem \ref{non-th:3.1} can be proved by some {\it a
priori} estimates globally in time.

For any given $T\in(0,+\infty)$, let $(\rho,u,\theta)$ be the
solution to (\ref{non-1.2})-(\ref{non-1.4}) as in Theorem
\ref{non-th:3.1}. Then we have the following basic energy estimate.

\begin{Lemma} \label{non-le:3.1} Under the conditions of Theorem \ref{non-th:3.1}, it holds for any $0\le t\le T$
\beqno &&\int_I \rho\left(1+e_c(\rho)+\theta^{1+r}+u^2\right)(t)\le
c. \eeqno
\end{Lemma}
\noindent{\it Proof.}  Integrating $(\ref{non-1.2})_1$ and
$(\ref{non-1.2})_3$ over $I\times [0,t]$, and using (\ref{non-1.4})
, $(A_2)$ and $(A_4)$, we complete the proof of Lemma
\ref{non-le:3.1}. \qed
\begin{Lemma} \label{non-le:3.2} Under the conditions of Theorem \ref{non-th:3.1}, it holds for any $(x, t)\in Q_T$
\beqno \begin{cases}0<\rho(x, t)\le c,\\ \theta(x,t)>0.\end{cases}
\eeqno
\end{Lemma}
\noindent {\it Proof.} The proof of the upper bound of $\rho$ relies
on constant viscosity (i.e. $\mu=const.$). It is similar to
\cite{Zhang Jianwen}.

Denote \be\label{non-3.1} w(x,t)=\int_0^t(u_x-P-\rho
u^2)+\int_0^x\rho_0u_0. \ee Differentiating (\ref{non-3.1}) with
respect to $x$, and using $(\ref{non-1.2})_2$, we have
$$
w_x=\rho u.
$$
This together with Lemma \ref{non-le:3.1} and the Cauchy inequality
gives
$$
\int_I|w_x|\le c.
$$
It follows from (\ref{non-3.1}), (\ref{non-1.1}), ($A_2$), ($A_3$),
($A_4$), (\ref{non-1.4}), and Lemma \ref{non-le:3.1} that
$$
\left|\int_Iw\right|\le c.
$$
This gives for any $(x,t)\in Q_T$ \beqno |w(x,t)|&\le &
\left|w(x,t)-\int_Iw\right|+\left|\int_Iw\right| \\ [4mm] & \le &
\left|\int_I\int_y^xw_\xi(\xi,t)d\xi dy\right|+c
\\ [4mm] &\le&\int_I|w_x|+c\le c, \eeqno which implies \be\label{non-3.2}
\|w\|_{L^\infty(Q_T)}\le c. \ee For any $(x,t)\in Q_T$, let $X(s;
x,t)$ satisfy
 \beq\label{non-3.3}
\begin{cases}
\frac{d X(s; x,t)}{d s}=u\left(X(s; x,t),s\right),\ 0\le s<t,\\
X(t; x,t)=x.
\end{cases} \eeq
Denote
$$
F(x,t)=\exp\left\{w(x,t)\right\}.
$$
It is easy to verify \beq\label{non-3.4}\nonumber \frac{d(\rho
F)\left(X(s; x,t), s\right)}{ds}&=&F\left(\rho_s+\frac{\partial
\rho}{\partial X}u+\rho \frac{\partial w}{\partial X}u+\rho
w_s\right)\\&=&-\rho PF. \eeq Multiplying (\ref{non-3.4}) by
$\exp\left(\int_0^sP\right)$, we have
$$
\frac{d}{ds}\left\{\rho F\exp\left(\int_0^sP\right)\right\}=0.
$$
Integrating it over $(0,t)$, we have \be\label{non-r3.5}
\rho(x,t)=\frac{F(X(0;
x,t),0)}{F(x,t)}\rho_0^\delta\exp\left(-\int_0^tP\right), \ee which
implies
$$
\rho(x,t)>0,
$$
for any $(x,t)\in Q_T$.

By (\ref{non-3.2}), (\ref{non-r3.5}) and $P\ge0$, we get the upper
bound of $\rho$. The lower bound of $\theta$ can be got by
(\ref{non-3.6}) and the maximum principle for parabolic
equations.\qed
\begin{Lemma}\label{non-le:3.3}
Under the conditions of Theorem \ref{non-th:3.1}, it holds for any
given $\alpha\in(0,1)$
$$
\int_{Q_T}\left(\frac{u_x^2}{\theta^\alpha}+\frac{(1+\theta^q)\theta_x^2}{\theta^{1+\alpha}}\right)\le
c,
$$
where $c$ may depend on $\alpha$.
\end{Lemma}
\begin{Remark}
$\alpha$ was usually taken as $1$ when the basic energy inequality
was done (see \cite{Bresch-Desjardins} and references therein). This
depends on $\rho_0\log\theta_0\in L^1$ which can not be got under
the assumptions of Theorem \ref{non-rth:1.1} and Theorem
\ref{non-th:1.1}, since $\theta_0$ may vanish.
\end{Remark}
\noindent {\it Proof.} From (\ref{non-1.1}) and (\ref{non-1.2}), we
get \be\label{non-3.5} \rho e_\theta\theta_t+\rho u e_\theta
\theta_x+\theta P_\theta
u_x=u_x^2+\left(\kappa(\theta)\theta_x\right)_x. \ee Substituting
$e=Q(\theta)+e_c(\rho)$ and $P=\rho Q(\theta)+P_c(\rho)$ into
(\ref{non-3.5}), we get \be\label{non-3.6} \rho
Q'(\theta)\theta_t+\rho u Q'(\theta)\theta_x+\rho\theta
Q'(\theta)u_x=u_x^2+\left(\kappa(\theta)\theta_x\right)_x, \ee or
\be\label{non-3.7} (\rho Q)_t+(\rho u Q)_x+\rho\theta
Q'(\theta)u_x=u_x^2+\left(\kappa(\theta)\theta_x\right)_x. \ee
Multiplying (\ref{non-3.6}) by $\theta^{-\alpha}$, integrating the
resulting equation over $Q_T$, and using integration by parts, we
have
\beq\label{non-3.8}\nonumber&&\int_{Q_T}\left(\frac{u_x^2}{\theta^\alpha}+\frac{\alpha\kappa(\theta)\theta_x^2}{\theta^{1+\alpha}}\right)\\&=&\nonumber
\int_I\rho\int_0^\theta\frac{Q'(\xi)}{\xi^\alpha}-\int_I\rho_0\int_0^{\theta_0}\frac{Q'(\xi)}{\xi^\alpha}+\int_{Q_T}\rho\theta^{1-\alpha}Q'(\theta)u_x\\&\le&\nonumber
c\int_I\int_0^\theta\frac{1+\xi^r}{\xi^\alpha}+c\int_I\rho_0\int_0^{\theta_0}\xi^{r-\alpha}+c\int_{Q_T}\rho\theta^{1-\alpha}(1+\theta^r)|u_x|\\&\le&\nonumber
c\int_I\rho(1+\theta^{1+r})+c+\frac{1}{2}\int_{Q_T}\frac{u_x^2}{\theta^\alpha}+c\int_{Q_T}\rho^2\theta^{2-\alpha+2r}\\&\le&c+\frac{1}{2}\int_{Q_T}\frac{u_x^2}{\theta^\alpha}
+c\int_0^T\max\limits_{x\in I}\theta^{1+r-\alpha}, \eeq where we
have used ($A_4$), the Cauchy inequality, Lemma \ref{non-le:3.1} and
Lemma \ref{non-le:3.2}. Now we estimate the last term of
(\ref{non-3.8}) as follows:
\beq\label{non-3.9} \nonumber
c\int_0^T\max\limits_{x\in I}\theta^{1+r-\alpha}&\le&
c+\int_0^T\|\theta^{r-\alpha}\theta_x\|_{L^2}\\&\le&\nonumber
c+c\int_0^T\left(\int_I\frac{\theta_x^2\theta^{2r-\alpha+1}}{\theta^{1+\alpha}}\right)^\frac{1}{2}\\&\le&c+
\frac{1}{2}\int_{Q_T}\frac{\alpha\kappa(\theta)\theta_x^2}{\theta^{1+\alpha}},
\eeq where we have used Corollary \ref{non-cor:2.1}, Lemma
\ref{non-le:3.1}, ($A_5$) and the Cauchy inequality. By
(\ref{non-3.8}), (\ref{non-3.9}) and ($A_5$), we complete the proof.
\qed

\begin{cor}\label{non-rcor:3.1}Under the conditions of Theorem \ref{non-th:3.1}, it holds
\beqno  \int_0^T\|\theta\|_{L^\infty}^{q-\alpha+1}\le c.\eeqno
\end{cor}
\noindent {\it Proof.} By Corollary \ref{non-cor:2.1} and Lemma
\ref{non-le:3.1}, we have \beqno
\int_0^T\|\theta\|_{L^\infty}^{q-\alpha+1}&=&\int_0^T\|\theta^\frac{q-\alpha+1}{2}\|_{L^\infty}^2\\&\le&
c\int_0^T\int_I\left(\theta^{\frac{q-\alpha-1}{2}}\theta_x\right)^2+c\\&=&c\int_0^T\int_I\theta^{q-\alpha-1}\theta_x^2+c\\&\le&
c. \eeqno \qed

\begin{Lemma}\label{non-le:3.5}Under the conditions of Theorem \ref{non-th:3.1}, it holds
$$\int_{Q_T}u_x^2\le c.$$
\end{Lemma}
\noindent {\it Proof.} From (\ref{non-1.2})$_1$ and
(\ref{non-1.2})$_2$, we get \be\label{non-3.10} \rho u_t+\rho
uu_x+P_x=u_{xx}. \ee Multiplying (\ref{non-3.10}) by $u$,
integrating it over $I$, and using integration by parts, we have
\beqno \frac{1}{2}\frac{d}{dt}\int_I\rho
u^2+\int_Iu_x^2&=&\int_IPu_x\\&\le&\frac{1}{2}\int_Iu_x^2+c\int_I\rho^2Q^2+c\int_IP_c^2\\&\le&\frac{1}{2}\int_Iu_x^2+c\int_I\theta^{2+2r}+c,\eeqno
where we have used the Cauchy inequality, ($A_2$), ($A_3$), ($A_4$)
and Lemma \ref{non-le:3.2}. This implies \beqno
\frac{d}{dt}\int_I\rho u^2+\int_Iu_x^2\le c\sup\limits_{x\in
I}\theta^{q-\alpha+1}+c. \eeqno Integrating it over $(0, t)$, and
using Corollary \ref{non-rcor:3.1}, we complete the proof of Lemma
3.4. \qed

\begin{Lemma} \label{non-le:3.6} Under the conditions of Theorem \ref{non-th:3.1}, it holds for any $0\le t\le T$
$$ \int_I(u_x^2+\rho\theta^{q+2+r})+\int_{Q_T}\left(\rho u_t^2+(1+\theta^q)^2\theta_x^2\right)\le c. $$
\end{Lemma}
\noindent {\it Proof.}
 Multiplying (\ref{non-3.10}) by $u_t$, integrating it over $I$, and using integration by parts, Lemma \ref{non-le:2.2}, Lemma \ref{non-le:3.2} and the Cauchy inequality, we have
 \beqno
\int_I\rho
u_t^2+\frac{1}{2}\frac{d}{dt}\int_Iu_x^2&=&\frac{d}{dt}\int_IP
u_x-\int_I\rho uu_xu_t-\int_IP_tu_x\\&\le&\frac{1}{2}\int_I\rho
u_t^2+\frac{1}{2}\int_I\rho u^2u_x^2+\frac{d}{dt}\int_IP
u_x-\int_IP_tu_x\\&\le&\frac{1}{2}\int_I\rho
u_t^2+c\left(\int_Iu_x^2\right)^2+\frac{d}{dt}\int_IP
u_x-\int_IP_t(u_x-P)-\frac{1}{2}\frac{d}{dt}\int_IP^2,
 \eeqno
which implies \beq\label{non-3.11} \int_I\rho
u_t^2+\frac{d}{dt}\int_Iu_x^2\le
c\left(\int_Iu_x^2\right)^2+2\frac{d}{dt}\int_IP
u_x-\frac{d}{dt}\int_IP^2-2\int_IP_t(u_x-P).\eeq We are going to
estimate the last term of the right side of (\ref{non-3.11}). Using
($A_2$), (\ref{non-3.7}), (\ref{non-1.2})$_1$ and integration by
parts, we have \beqno -2\int_IP_t(u_x-P)&=&-2\int_I(\rho
Q)_t(u_x-P)-2\int_I(P_c)_t(u_x-P)\\&=&-2\int_I\left[(\kappa
\theta_x)_x+u_x^2-(\rho uQ)_x-\rho\theta
Q'(\theta)u_x\right](u_x-P)\\&&+2\int_IP_c^\prime(\rho)(\rho_xu+\rho
u_x)(u_x-P)\\&=&2\int_I\kappa\theta_x(u_{xx}-P_x)-2\int_Iu_x^2(u_x-P)-2\int_I\rho
uQ(u_{xx}-P_x)\\&&+2\int_I\rho\theta
Q^\prime(\theta)u_x(u_x-P)-2\int_IP_cu(u_{xx}-P_x)-2\int_IP_cu_x(u_x-P)\\&&+2\int_I\rho
P_c^\prime(\rho)u_x(u_x-P). \eeqno  This, together with
(\ref{non-3.10}), ($A_2$), ($A_4$), Lemma \ref{non-le:2.2}, Lemma
\ref{non-le:3.2}, the Cauchy inequality, and
$W^{1,1}(I)\hookrightarrow L^\infty(I)$, gives
\beq\label{non-3.12}\nonumber
-2\int_IP_t(u_x-P)&\le&2\int_I\kappa\theta_x(\rho u_t+\rho
uu_x)+2\|u_x-P\|_{L^\infty}\int_Iu_x^2-2\int_I\rho u Q(\rho u_t+\rho
uu_x)\\&&\nonumber+c\sup\limits_{x\in
I}(1+\theta^{1+r})\int_Iu_x^2+c\sup\limits_{x\in
I}(1+\theta^{1+r})\int_I\rho Q^2+c\sup\limits_{x\in
I}\theta^{1+r}\\&&\nonumber-2\int_IP_cu(\rho u_t+\rho
uu_x)+c\int_Iu_x^2+c\int_I\rho
Q^2+c\\&\le&\nonumber\frac{1}{4}\int_I\rho
u_t^2+c\int_I\kappa^2\theta_x^2+c\left(\int_Iu_x^2\right)^2+c\left(\|u_x-P\|_{L^1}+\|\rho
u_t+\rho
uu_x\|_{L^1}\right)\int_Iu_x^2\\&&\nonumber+c\int_Iu_x^2\int_I\rho
Q^2+c\sup\limits_{x\in I}(1+\theta^{1+r})\int_I(u_x^2+\rho
Q^2)+c\sup\limits_{x\in
I}\theta^{1+r}+c\\&\le&\nonumber\frac{1}{4}\int_I\rho
u_t^2+c\int_I\kappa^2\theta_x^2+c\left(\int_Iu_x^2\right)^2+\frac{1}{4}\int_I\rho
u_t^2+c\int_Iu_x^2\int_I\rho Q^2\\&&+c\sup\limits_{x\in
I}(1+\theta^{1+r})\int_I(u_x^2+\rho Q^2)+c\sup\limits_{x\in
I}\theta^{1+r}+c. \eeq
Substituting (\ref{non-3.12}) into (\ref{non-3.11}), we have
\beq\label{non-3.13}\nonumber&&\frac{1}{2} \int_I\rho
u_t^2+\frac{d}{dt}\int_Iu_x^2\\&\le&\nonumber
c\left(\int_Iu_x^2\right)^2+2\frac{d}{dt}\int_IP
u_x-\frac{d}{dt}\int_IP^2+c\int_I\kappa^2\theta_x^2+c\int_Iu_x^2\int_I\rho
Q^2\\&&+c\sup\limits_{x\in I}(1+\theta^{1+r})\int_I(u_x^2+\rho
Q^2)+c\sup\limits_{x\in I}\theta^{1+r}+c.\eeq Integrating
(\ref{non-3.13}) over $(0, t)$, and using ($A_2$)-($A_4$), Lemma
\ref{non-le:3.2}, Corollary \ref{non-rcor:3.1}, Lemma
\ref{non-le:3.5} and the Cauchy inequality, we have
\beqno&&\frac{1}{2} \int_0^t\int_I\rho u_t^2+\int_Iu_x^2\\&\le&
c\int_0^t\left(\int_Iu_x^2\right)^2+2\int_I(\rho Q+P_c)
u_x+c\int_0^t\int_I\kappa^2\theta_x^2+c\int_0^t\int_Iu_x^2\int_I\rho
\theta^{2+2r}\\&&+c\int_0^t\sup\limits_{x\in
I}\theta^{1+r}\int_Iu_x^2+c\int_0^t\sup\limits_{x\in
I}(1+\theta^{1+r})\int_I\rho \theta^{2+2r}+c\\&\le&
c\int_0^t\left(\int_Iu_x^2\right)^2+\frac{1}{2}\int_Iu_x^2+c\int_I\rho\theta^{2+2r}+c\int_0^t\int_I\kappa^2\theta_x^2+c\int_0^t\int_Iu_x^2\int_I\rho
\theta^{2+2r}\\&&+c\int_0^t\sup\limits_{x\in
I}\theta^{1+r}\int_Iu_x^2+c\int_0^t\sup\limits_{x\in
I}(1+\theta^{1+r})\int_I\rho \theta^{2+2r}+c.\eeqno The second term
of the right side can be absorbed by the left. After that, we have
\beq\label{non-3.14} \nonumber&&\int_0^t\int_I\rho
u_t^2+\int_Iu_x^2\\&\le&\nonumber
c\int_0^t\left(\int_Iu_x^2\right)^2+c\int_I\rho\theta^{q+2+r}+c\int_0^t\int_I\kappa^2\theta_x^2+c\int_0^t\int_Iu_x^2\int_I\rho
\theta^{2+2r}\\&&+c\int_0^t\sup\limits_{x\in
I}\theta^{1+r}\int_Iu_x^2+c\int_0^t\sup\limits_{x\in
I}(1+\theta^{1+r})\int_I\rho \theta^{2+2r}+c. \eeq Here, we have
used Lemma \ref{non-le:3.2} and the Young inequality on the second
term of the right side. Note that the terms about $\theta$ in
(\ref{non-3.14}) need to be handled. To do this, we make use of
($\ref{non-3.6}$).

 Multiplying
(\ref{non-3.6}) by $\int_0^\theta\kappa(\xi)d\xi$, integrating it
over $I$, and using integration by parts, ($A_4$) and ($A_5$), we
have \beq\label{non-3.15}
\nonumber&&\frac{d}{dt}\int_I\rho\left[\int_0^\theta
Q^\prime(\eta)\int_0^\eta\kappa(\xi)d\xi
d\eta\right]+\int_I\kappa^2\theta_x^2\\&&\nonumber=\int_Iu_x^2\int_0^\theta\kappa(\xi
)d\xi-\int_I\rho\theta
Q^\prime(\theta)u_x\int_0^\theta\kappa(\xi)d\xi\\&&\le
c\|(1+\theta^q)\theta\|_{L^\infty}\int_Iu_x^2+c\|(1+\theta^q)\theta\|_{L^\infty}\int_I\rho(1+\theta^{1+r})|u_x|.
\eeq By Corollary \ref{non-cor:2.1} and ($A_5$), we get
\be\label{non-3.16}\|(1+\theta^q)\theta\|_{L^\infty}\le
c\|\kappa\theta_x\|_{L^2}+c.\ee Substituting (\ref{non-3.16}) into
(\ref{non-3.15}), and using the H\"older inequality, the Cauchy
inequality and Lemma \ref{non-le:3.2}, we get \beqno
&&\frac{d}{dt}\int_I\rho\left[\int_0^\theta
Q^\prime(\eta)\int_0^\eta\kappa(\xi)d\xi
d\eta\right]+\int_I\kappa^2\theta_x^2\\&\le&c\|\kappa\theta_x\|_{L^2}\int_Iu_x^2+c\int_Iu_x^2+c\|\kappa\theta_x\|_{L^2}\int_I\rho(1+\theta^{1+r})|u_x|
+c\int_I\rho(1+\theta^{1+r})|u_x|\\&\le&c\|\kappa\theta_x\|_{L^2}\left(\int_Iu_x^2+\|\rho(1+\theta^{1+r})\|_{L^2}\|u_x\|_{L^2}\right)+c\int_Iu_x^2
+c\int_I\rho(1+\theta^{2+2r})\\&\le&\frac{1}{2}\int_I\kappa^2\theta_x^2+c\left(\int_Iu_x^2\right)^2+c\int_I\rho\theta^{2+2r}\int_Iu_x^2+c\int_I\rho\theta^{2+2r}+c,
\eeqno which implies \beqno
&&\frac{d}{dt}\int_I\rho\left[\int_0^\theta
Q^\prime(\eta)\int_0^\eta\kappa(\xi)d\xi
d\eta\right]+\frac{1}{2}\int_I\kappa^2\theta_x^2\\&\le&c\left(\int_Iu_x^2\right)^2+c\int_I\rho\theta^{2+2r}\int_Iu_x^2+c\int_I\rho\theta^{2+2r}+c.
\eeqno Integrating it over $(0, t)$, and using ($A_4$), ($A_5$),
Lemma \ref{non-le:3.1} and Corollary \ref{non-rcor:3.1}, we get
\be\label{non-3.17}
\int_I\rho\theta^{q+2+r}+\int_0^t\int_I\kappa^2\theta_x^2\le
c\int_0^t\left(\int_Iu_x^2\right)^2+c\int_0^t\left(\int_I\rho\theta^{2+2r}\int_Iu_x^2\right)+c.
\ee By (\ref{non-3.14}), (\ref{non-3.17}), Corollary
\ref{non-rcor:3.1}, Lemma \ref{non-le:3.5}, and the Gronwall
inequality, we complete the proof. \qed

\begin{Lemma} \label{non-le:3.7} Under the conditions of Theorem \ref{non-th:3.1}, it holds for any $0\le t\le T$
$$
\int_I(\rho_x^2+\rho_t^2)+\int_{Q_T}u_{xx}^2\le c.
$$
\end{Lemma}
\noindent {\it Proof.} Differentiating $(\ref{non-1.2})_1$ with
respect to $x$, we have \be\label{non-3.18}
\rho_{xt}+\rho_{xx}u+2\rho_xu_x+\rho u_{xx}=0. \ee Multiplying
(\ref{non-3.18}) by $2\rho_x$, integrating it over $I$ and using
integration by parts, we have \beq\label{non-3.19}
\nonumber\frac{d}{dt}\int_I\rho_x^2&=&
-3\int_I\rho_x^2u_x-2\int_I\rho\rho_xu_{xx}\\&=&\nonumber-3\int_I\rho_x^2(u_x-P)-3\int_I\rho_x^2P-2\int_I\rho\rho_xu_{xx}\\&\le&\nonumber
c\left(\|u_x-P\|_{L^2}+\|u_{xx}-P_x\|_{L^2}\right)
\int_I\rho_x^2+c\int_Iu_{xx}^2+c\int_I\rho_x^2\\&\le&
c\left(1+\|\sqrt{\rho}u_t\|_{L^2}\right)\int_I\rho_x^2+c\int_Iu_{xx}^2,
\eeq where we have used (\ref{non-3.10}), the Sobolev inequality,
($A_2$)-($A_4$), the Cauchy inequality, Lemma \ref{non-le:2.2},
Lemma \ref{non-le:3.2} and Lemma \ref{non-le:3.6}.

 It follows from (\ref{non-3.10}), Lemma \ref{non-le:2.2}, Lemma \ref{non-le:3.2}, ($A_2$)-($A_4$), Lemma \ref{non-le:3.6} and the Cauchy inequality that \beq\label{non-3.20}\nonumber \int_Iu_{xx}^2&\le&\nonumber c\int_I\rho
u_t^2+c\left(\int_Iu_x^2\right)^2 +c\int_I\rho_x^2Q^2
+c\int_I\rho^2\left[Q^\prime(\theta)\right]^2\theta_x^2+c\int_I\rho_x^2+c\\&\le&\nonumber
c\int_I\rho u_t^2+c\sup\limits_{x\in
I}(1+\theta^{2+2r})\int_I\rho_x^2+c\int_I(1+\theta^q)^2\theta_x^2+c\\&\le&c\int_I\rho
u_t^2+c\sup\limits_{x\in
I}(1+\theta^{q-\alpha+1})\int_I\rho_x^2+c\int_I(1+\theta^q)^2\theta_x^2+c.
\eeq Substituting (\ref{non-3.20}) into (\ref{non-3.19}), and using
the Gronwall inequality, Corollary \ref{non-rcor:3.1} and Lemma
\ref{non-le:3.6}, we get \be\label{non-3.21} \int_I\rho_x^2\le c.
\ee By (\ref{non-3.20}), (\ref{non-3.21}), Corollary
\ref{non-rcor:3.1} and Lemma \ref{non-le:3.6}, we have \beqno
\int_{Q_T}u_{xx}^2\le c. \eeqno It follows from (\ref{non-1.2})$_1$,
(\ref{non-3.21}), Lemma \ref{non-le:2.2}, Lemma \ref{non-le:3.2} and
Lemma \ref{non-le:3.6} that
$$
\int_I\rho_t^2\le c.
$$
The proof of the lemma is complete.
 \qed

\begin{Lemma} \label{non-le:3.8} Under the conditions of Theorem \ref{non-th:3.1}, it holds for any $0\le t\le T$
$$
\int_I\left(\rho
u_t^2+\theta_x^2\right)+\int_{Q_T}\left(u_{xt}^2+\rho\theta_t^2\right)\le
c.
$$
\end{Lemma}
\noindent{\it Proof.} Differentiating (\ref{non-3.10}) with respect
to $t$, we have \be\label{non-3.22} \rho
u_{tt}+\rho_tu_t+\rho_tuu_x+\rho u_tu_x+\rho uu_{xt}+P_{xt}=u_{xxt}.
\ee Multiplying (\ref{non-3.22}) by $u_t$, integrating the resulting
equation over $I$, we have \beqno
&&\frac{1}{2}\frac{d}{dt}\int_I\rho
u_t^2+\int_Iu_{xt}^2\\&=&-2\int_I\rho
uu_tu_{xt}-\int_I\rho_tuu_xu_t-\int_I\rho u_t^2u_x+\int_IP_tu_{tx}
\\&\le& 2\|\sqrt{\rho}u_t\|_{L^2}
\|\sqrt{\rho}u\|_{L^\infty}\|u_{xt}\|_{L^2}+\|u_t\|_{L^\infty}\|u\|_{L^\infty}\|\rho_t\|_{L^2}\|u_x\|_{L^2}
 +
\|u_x\|_{L^\infty}\int_I\rho
u_t^2\\&&+\|P_c^\prime(\rho)\|_{L^\infty}\|\rho_t\|_{L^2}\|u_{xt}\|_{L^2}+\|Q(\theta)\|_{L^\infty}\|\rho_t\|_{L^2}
\|u_{xt}\|_{L^2} +\|\rho
Q^\prime(\theta)\theta_t\|_{L^2}\|u_{xt}\|_{L^2}\\&\le&\frac{1}{2}\int_Iu_{xt}^2+c\int_I\rho
u_t^2+c+c\int_Iu_{xx}^2\int_I\rho u_t^2+c\sup\limits_{x\in
I}\theta^{2+2r}+c\int_I\rho\left(1+\theta^{q+r}\right)\theta_t^2.
\eeqno Here, we have used $(\ref{non-1.2})_1$, integration by parts,
the H\"older inequality, the Cauchy inequality, the Sobolev
inequality, ($A_2$)-($A_4$), Lemma \ref{non-le:2.2}, Lemma
\ref{non-le:3.2}, Lemma \ref{non-le:3.6} and Lemma \ref{non-le:3.7}.

The first term of the right side can be absorbed by the left. This
implies \beq\label{non-3.23} &&\nonumber\frac{d}{dt}\int_I\rho
u_t^2+\int_Iu_{xt}^2\\&\le&c\int_I\rho
u_t^2+c+c\int_Iu_{xx}^2\int_I\rho u_t^2+c\sup\limits_{x\in
I}\theta^{2+2r}+c\int_I\rho\left(1+\theta^{q+r}\right)\theta_t^2.
\eeq Integrating (\ref{non-3.23}) over $(0, t)$, and using Corollary
\ref{non-rcor:3.1} and Lemma \ref{non-le:3.6}, we have
\beq\label{non-3.24} \int_I\rho u_t^2+\int_0^t\int_Iu_{xt}^2\le
\int_I\rho u_t^2(0)+c+c\int_0^t\int_Iu_{xx}^2\int_I\rho
u_t^2+c\int_0^t\int_I\rho\left(1+\theta^{q+r}\right)\theta_t^2. \eeq
Multiplying (\ref{non-3.10}) by $\frac{1}{\sqrt{\rho}}$, taking
$t\rightarrow0^+$ and using (\ref{non-1.8})$_1$, we have \beqno
|\sqrt{\rho}
u_t(x,0)|&\le&\frac{\left|u_{0xx}-P(\rho_0^\delta,\theta_0^\delta)_x\right|}{\sqrt{\rho_0^\delta}}+\sqrt{\rho_0^\delta}
|u_0u_{0x}|\\&\le&\frac{\left|u_{0xx}-P(\rho_0,\theta_0)_x\right|}{\sqrt{\rho_0^\delta}}
+\frac{|P(\rho_0,\theta_0)_x-P(\rho_0^\delta,\theta_0^\delta)_x|}{\sqrt{\rho_0^\delta}}+\sqrt{\rho_0^\delta}
|u_0u_{0x}|\\&\le&|g_1|+c\frac{\delta}{\sqrt{\rho_0^\delta}}(|\rho_{0x}|+|\theta_{0x}|)+c,
\eeqno which implies \be\label{non-rhou_t^2(0)} \int_I\rho
u_t^2(0)\le c. \ee Substituting (3.26) into (\ref{non-3.24}), we
have \beq\label{non-r3.24} \int_I\rho
u_t^2+\int_0^t\int_Iu_{xt}^2\le c+c\int_0^t\int_Iu_{xx}^2\int_I\rho
u_t^2+c\int_0^t\int_I\rho\left(1+\theta^{q+r}\right)\theta_t^2. \eeq
Multiplying (\ref{non-3.6}) by
$\left(\int_0^\theta\kappa(\xi)d\xi\right)_t$(i.e.
$\kappa(\theta)\theta_t$), integrating the resulting equation over
$I$, and using integration by parts, ($A_4$), ($A_5$), Lemma
\ref{non-le:2.2}, Lemma \ref{non-le:3.2}, Lemma \ref{non-le:3.6} and
the Cauchy inequality, we have for any $\varepsilon>0$ \beqno&&
\int_I\rho
Q^\prime(\theta)\kappa(\theta)\theta_t^2+\frac{1}{2}\frac{d}{dt}\int_I\kappa^2\theta_x^2\\&=&-\int_I\rho
uQ^\prime\kappa\theta_x\theta_t-\int_I\rho\theta Q^\prime\kappa
u_x\theta_t+\int_Iu_x^2\left(\int_0^\theta\kappa(\xi)d\xi\right)_t
\\ & \le & \frac{1}{2}\int_I\rho Q^\prime\kappa\theta_t^2+c\int_I\rho
u^2Q^\prime\kappa\theta_x^2+c\int_I\rho Q^2Q^\prime\kappa u_x^2 \\
&&
+\frac{d}{dt}\left(\int_Iu_x^2\int_0^\theta\kappa(\xi)d\xi\right)-2\int_Iu_xu_{xt}\int_0^\theta\kappa(\xi)d\xi\\&\le&\frac{1}{2}\int_I\rho
Q^\prime\kappa\theta_t^2+c\int_I(1+\theta^q)^2\theta_x^2+c\left(1+\int_Iu_{xx}^2\right)\int_I\rho(1+\theta^{q+r+2})
\\&&+\frac{d}{dt}\left(\int_Iu_x^2\int_0^\theta\kappa(\xi)d\xi
\right)+\varepsilon\int_Iu_{xt}^2+c_\varepsilon\sup\limits_{x\in
I}(1+\theta^q)^2\theta^2, \eeqno which combining Lemma
\ref{non-le:2.2}, Lemma \ref{non-le:3.2}, Lemma \ref{non-le:3.6}
implies \beqno && \int_I\rho
Q^\prime(\theta)\kappa(\theta)\theta_t^2+\frac{d}{dt}\int_I\kappa^2\theta_x^2\\&\le&c\int_I(1+\theta^q)^2\theta_x^2+c\int_Iu_{xx}^2+c
+\frac{d}{dt}\left(\int_Iu_x^2\int_0^\theta\kappa(\xi)d\xi \right)
+\varepsilon\int_Iu_{xt}^2+c_\varepsilon\sup\limits_{x\in
I}(1+\theta^q)^2\theta^2\\&\le&c\int_Iu_{xx}^2+\frac{d}{dt}\left(\int_Iu_x^2\int_0^\theta\kappa(\xi)d\xi
\right)+\varepsilon\int_Iu_{xt}^2
+c_\varepsilon\int_I(1+\theta^q)^2\theta_x^2+c_\varepsilon. \eeqno
Integrating it over $(0, t)$, and using ($A_4$) and ($A_5$), Lemma
\ref{non-le:2.2}, Lemma \ref{non-le:3.6}, Lemma \ref{non-le:3.7} and
the Cauchy inequality, we obtain
\beqno&&\int_0^t\int_I\rho\left(1+\theta^{q+r}\right)\theta_t^2+\int_I(1+\theta^q)^2\theta_x^2\\&\le&c\int_Iu_x^2\int_0^\theta\kappa(\xi)d\xi
+c\varepsilon\int_0^t\int_Iu_{xt}^2+c_\varepsilon\\&\le&c\sup\limits_{x\in
I}(1+\theta^q)\theta+c\varepsilon\int_0^t\int_Iu_{xt}^2+c_\varepsilon\\&\le&c\|(1+\theta^q)\theta_x\|_{L^2}+c\varepsilon\int_0^t\int_Iu_{xt}^2
+c_\varepsilon\\&\le&\frac{1}{2}\int_I(1+\theta^q)^2\theta_x^2+c\varepsilon\int_0^t\int_Iu_{xt}^2
+c_\varepsilon.\eeqno After the first term of the right side is
absorbed by the left, we get \be\label{non-3.25}
\int_0^t\int_I\rho\left(1+\theta^{q+r}\right)\theta_t^2+\int_I(1+\theta^q)^2\theta_x^2\le
c\varepsilon\int_0^t\int_Iu_{xt}^2 +c_\varepsilon. \ee Multiplying
(\ref{non-3.25}) by $2c$, adding the resulting inequality to
(\ref{non-3.24}), taking $\varepsilon=\frac{1}{4c^2}$, and using the
Gronwall inequality and Lemma \ref{non-le:3.7}, we complete the
proof of Lemma \ref{non-le:3.8}. \qed\\[1.5mm]

From Corollary \ref{non-cor:2.1}, Lemma \ref{non-le:3.1} and Lemma
\ref{non-le:3.8}, we get the following corollary immediately.
\begin{cor}\label{non-cor:3.1}Under the conditions of Theorem \ref{non-th:3.1}, it holds
$$\|\theta\|_{L^\infty(Q_T)}\le c.$$
\end{cor}
\begin{cor} \label{non-rcor:3.3} Under the conditions of Theorem \ref{non-th:3.1}, it holds for any $0\le t\le T$
$$
\|u\|_{W^{1,\infty}(Q_T)}+\int_I u_{xx}^2+\int_{Q_T}\theta_{xx}^2\le
c.
$$
\end{cor}
\noindent{\it Proof.} It follows from (\ref{non-3.20}), Lemma
\ref{non-le:3.7}, Lemma \ref{non-le:3.8} and Corollary
\ref{non-cor:3.1} that
$$
\int_Iu_{xx}^2\le c,
$$
which, combining Lemma \ref{non-le:2.2}, Lemma \ref{non-le:3.6} and
the Sobolev inequality, gives \be\label{non-r3.29}
\|u\|_{W^{1,\infty}(Q_T)}\le c. \ee
By (\ref{non-3.6}), Corollary \ref{non-cor:3.1}, ($A_4$), ($A_5$),
Lemma \ref{non-le:2.2}, Lemma \ref{non-le:3.2}, (\ref{non-r3.29}),
Lemma \ref{non-le:3.8}, the H\"older inequality, Sobolev inequality
and Cauchy inequality, we have \beqno \int_I\theta_{xx}^2&\le&
c\int_I\theta_x^4+c\int_Iu_x^4+\int_I\rho\theta_t^2+c\int_Iu^2\theta_x^2+c\int_I\theta^2u_x^2\\&\le&c\|\theta_x\theta_{xx}\|_{L^1}\int_I\theta_x^2
+\int_I\rho\theta_t^2+c\\&\le&c\|\theta_{xx}\|_{L^2}
+\int_I\rho\theta_t^2+c\\&\le&\frac{1}{2}\int_I\theta_{xx}^2
+\int_I\rho\theta_t^2+c.\eeqno After the first term of the right
side is absorbed by the left, we get \be\label{non-3.26}
\int_I\theta_{xx}^2\le \int_I\rho\theta_t^2+c. \ee Integrating
(\ref{non-3.26}) over $[0, T]$, and using Lemma \ref{non-le:3.8}, we
get
$$
\int_{Q_T}\theta_{xx}^2\le c.
$$
This proves Corollary 3.3. \qed
\begin{Lemma}\label{non-le:3.10} Under the conditions of Theorem \ref{non-th:3.1}, it holds for any $0\le t\le T$
$$
\|\rho\|_{W^{1,\infty}(Q_T)}+\|\rho_t\|_{L^\infty(Q_T)}+\int_I(\rho_{xx}^2+\rho_{xt}^2)+\int_{Q_T}(\rho_{tt}^2+u_{xxx}^2)\le
c.
$$
\end{Lemma}
\noindent{\it Proof.} Differentiating (\ref{non-3.18}) w.r.t. $x$,
we have \be\label{non-3.27}
\rho_{xxt}=-\rho_{xxx}u-3\rho_{xx}u_x-3\rho_xu_{xx}-\rho u_{xxx}.
\ee Multiplying (\ref{non-3.27}) by $2\rho_{xx}$, integrating it
over $I$, and using integration by parts and the H\"older
inequality, we have \beqno \frac{d}{dt}\int_I\rho_{xx}^2&=&
-5\int_I\rho_{xx}^2u_x-6\int_I\rho_x\rho_{xx}u_{xx}
-2\int_I\rho\rho_{xx}u_{xxx}\\&\le&5\|u_x\|_{L^\infty}\int_I\rho_{xx}^2+6\|\rho_x\|_{L^\infty}\|\rho_{xx}\|_{L^2}
\|u_{xx}\|_{L^2}+2\|\rho\|_{L^\infty}\|\rho_{xx}\|_{L^2}\|u_{xxx}\|_{L^2}.
\eeqno
By the Sobolev inequality, Cauchy inequality, Lemma
\ref{non-le:3.2}, Lemma \ref{non-le:3.7}, and Corollary
\ref{non-rcor:3.3}, we have \beq\label{non-3.28}
\frac{d}{dt}\int_I\rho_{xx}^2\le
c\int_I\rho_{xx}^2+c\int_Iu_{xxx}^2+c. \eeq
 The next step is to
estimate the term $\int_Iu_{xxx}^2$. Differentiating (\ref{non-3.10}) with respect to $x$, we have
\beq\label{non-3.29} u_{xxx}=\rho_xu_t+\rho u_{xt}+\rho_x uu_x+\rho
u_x^2+\rho uu_{xx}+(P_c)_{xx}+(\rho Q)_{xx}.\eeq By ($A_3$),
($A_4$), Lemma \ref{non-le:2.2}, Lemma \ref{non-le:3.2}, Lemma
\ref{non-le:3.7}, Corollary \ref{non-cor:3.1}, Corollary
\ref{non-rcor:3.3} and the Sobolev inequality, we get
\beq\label{non-3.30} \nonumber\int_Iu_{xxx}^2&\le&
c\int_I\rho^2u_{xt}^2+c\int_I\rho_x^2u_t^2
+c\int_I\rho_{xx}^2+c\int_I\theta_{xx}^2+c\\&\le&c\int_Iu_{xt}^2+c\int_I\rho_{xx}^2+c\int_I\theta_{xx}^2+c.
\eeq Substituting (\ref{non-3.30}) into (\ref{non-3.28}), and using
the Gronwall inequality, Lemma \ref{non-le:3.8} and Corollary
\ref{non-rcor:3.3}, we get \be\label{non-3.31} \int_I\rho_{xx}^2\le
c. \ee By (\ref{non-3.31}), Lemma \ref{non-le:3.2}, Lemma
\ref{non-le:3.7} and the Sobolev inequality, we have
\be\label{non-r3.36} \|\rho\|_{W^{1,\infty}(Q_T)}\le c. \ee
By (\ref{non-3.30}), (\ref{non-3.31}), Lemma \ref{non-le:3.8} and
Corollary \ref{non-rcor:3.3}, we get
$$
\int_{Q_T}u_{xxx}^2\le c.
$$
The estimates of $\rho_{xt}$ and $\rho_{tt}$ can be obtained
directly by (\ref{non-3.18}), (\ref{non-1.2})$_1$, (\ref{non-3.31}),
(\ref{non-r3.36}), Lemma \ref{non-le:2.2}, Lemma \ref{non-le:3.2},
Lemma \ref{non-le:3.7}, Lemma \ref{non-le:3.8}, and Corollary
\ref{non-rcor:3.3}. The proof of Lemma \ref{non-le:3.10} is
complete. \qed
\begin{Lemma}\label{non-rle:3.11}Under the conditions of Theorem \ref{non-th:3.1}, it holds for any $0\le t\le T$
$$
\int_I\rho\theta_t^2+\int_{Q_T}\left|(\kappa\theta_x)_t\right|^2\le
c.
$$
\end{Lemma}
 Differentiating
(\ref{non-3.6}) w.r.t. $t$, we have \beq\label{non-3.36} \rho
Q^\prime\theta_{tt}+\rho
Q^{\prime\prime}\theta_t^2+\rho_tQ^\prime\theta_t+(\rho
uQ^\prime\theta_x)_t+(\rho\theta Q^\prime
u_x)_t=2u_xu_{xt}+(\kappa\theta_x)_{xt}. \eeq Multiplying
(\ref{non-3.36}) by $(\int_0^\theta\kappa(\xi)d\xi)_t$ (i.e.
$\kappa(\theta)\theta_t$),
 integrating it over $I$, and using integration by parts, (\ref{non-1.2})$_1$, ($A_4$), ($A_5$), Corollary \ref{non-cor:3.1}, Lemma \ref{non-le:3.2},
  Corollary \ref{non-rcor:3.3} and the H\"older inequality, we have
\beqno &&\frac{1}{2}\frac{d}{dt}\int_I\rho
 Q^\prime\kappa\theta_t^2+\int_I\left|(\kappa\theta_x)_t\right|^2\\&=&-\frac{1}{2}\int_I\rho_tQ^\prime\kappa\theta_t^2-\frac{1}{2}\int_I\rho
 Q^{\prime\prime}\theta_t^3\kappa+\frac{1}{2}\int_I\rho
 Q^\prime\kappa^\prime\theta_t^3-\int_I(\rho u
 Q^\prime\theta_x)_t\kappa\theta_t \\&&-\int_I(\rho\theta Q^\prime
 u_x)_t\kappa\theta_t+2\int_Iu_xu_{xt}\kappa\theta_t\\&\le&\frac{1}{2}\int_I(\rho
 u)_xQ^\prime\kappa\theta_t^2+c\|\kappa\theta_t\|_{L^\infty}\int_I\rho\theta_t^2-\int_I\rho
 u Q^\prime(\kappa\theta_x)_t\theta_t-\int_I\rho
 u(Q^{\prime\prime}\kappa-Q^\prime\kappa^\prime)\theta_t^2\theta_x\\&&-\int_I(\rho
 u)_tQ^\prime\theta_x\kappa\theta_t+c\int_Iu_{xt}^2+c\int_I\rho\theta_t^2-\int_I\rho_t\theta
 Q^\prime
 u_x\kappa\theta_t+c\|\kappa\theta_t\|_{L^\infty}\|u_{xt}\|_{L^2}\|u_x\|_{L^2}.
\eeqno This, combining integration by parts, ($A_4$), ($A_5$), Lemma
\ref{non-le:2.1}, Lemma \ref{non-le:2.2}, Corollary
\ref{non-cor:3.1}, Corollary \ref{non-rcor:3.3}, Lemma
\ref{non-le:3.8}, Lemma \ref{non-le:3.10} and the Cauchy inequality,
gives \beqno &&\frac{1}{2}\frac{d}{dt}\int_I\rho
 Q^\prime\kappa\theta_t^2+\int_I\left|(\kappa\theta_x)_t\right|^2 \\&\le&-\frac{1}{2}\int_I\rho
 u Q^{\prime\prime}\theta_x\kappa\theta_t^2-\frac{1}{2}\int_I\rho u
 Q^\prime\kappa^\prime\theta_x\theta_t^2-\int_I\rho uQ^\prime\kappa\theta_t\theta_{xt}+c\|\kappa\theta_t\|_{L^\infty}\int_I\rho\theta_t^2
 \\&&+c\|\sqrt{\rho}\theta_t\|_{L^2}\|(\kappa\theta_x)_t\|_{L^2}
 +c\|\theta_x\|_{L^\infty}\int_I\rho\theta_t^2+c\|\kappa\theta_t\|_{L^\infty}+c\int_Iu_{xt}^2+c\int_I\rho\theta_t^2\\&&+c\|\kappa\theta_t\|_{L^\infty}\|u_{xt}\|_{L^2}
 \\&\le&c\|\theta_{xx}\|_{L^2}\int_I\rho\theta_t^2+c\|\kappa\theta_t\|_{L^\infty}\int_I\rho\theta_t^2+c\|\sqrt{\rho}\theta_t\|_{L^2}\|(\kappa\theta_x)_t\|_{L^2}
+c\|\kappa\theta_t\|_{L^\infty}
\\&&+c\int_Iu_{xt}^2+c\int_I\rho\theta_t^2+c\|\kappa\theta_t\|_{L^\infty}\|u_{xt}\|_{L^2}\\&\le&c\|\theta_{xx}\|_{L^2}\int_I\rho\theta_t^2
+c\left(\|(\kappa\theta_t)_x\|_{L^2}+\int_I\rho\kappa|\theta_t|\right)\int_I\rho\theta_t^2+c\|\sqrt{\rho}\theta_t\|_{L^2}\|(\kappa\theta_x)_t\|_{L^2}
\\&&+c\|(\kappa\theta_t)_x\|_{L^2}+c\int_I\rho\kappa|\theta_t|
+c\int_Iu_{xt}^2+c\int_I\rho\theta_t^2+c\left(\|(\kappa\theta_t)_x\|_{L^2}+\int_I\rho\kappa|\theta_t|\right)\|u_{xt}\|_{L^2},
\eeqno which together with the Cauchy inequality, Lemma
\ref{non-le:3.2}, ($A_5$) and Corollary \ref{non-cor:3.1} gives
\beqno \frac{1}{2}\frac{d}{dt}\int_I\rho
 Q^\prime\kappa\theta_t^2+\int_I\left|(\kappa\theta_x)_t\right|^2\le\frac{1}{2}\int_I\left|(\kappa\theta_x)_t\right|^2+c\int_I\theta_{xx}^2
 +c\left(\int_I\rho\theta_t^2\right)^2+c\int_Iu_{xt}^2+c,
\eeqno where we have used $(\kappa\theta_t)_x=(\kappa\theta_x)_t$.
This gives \beqno \frac{d}{dt}\int_I\rho
 Q^\prime\kappa\theta_t^2+\int_I\left|(\kappa\theta_x)_t\right|^2\le c\int_I\theta_{xx}^2
 +c\left(\int_I\rho\theta_t^2\right)^2+c\int_Iu_{xt}^2+c.
\eeqno Integrating it over $(0, t)$, and using ($A_4$), ($A_5$),
Lemma \ref{non-le:3.8} and Corollary \ref{non-rcor:3.3}, we obtain
\be\label{non-3.37} \int_I\rho
 \theta_t^2+\int_0^t\int_I\left|(\kappa\theta_x)_t\right|^2\le c\int_I\rho\theta_t^2(0)+c\int_0^t\left(\int_I\rho\theta_t^2\right)^2+c.
\ee Multiplying (\ref{non-3.6}) by $\D\frac{1}{
Q^\prime(\theta)\sqrt{\rho}}$, taking $t\rightarrow0^+$, and using
(\ref{non-1.8})$_2$, we have \beqno
|\sqrt{\rho}\theta_t(x,0)|&\le&\frac{\left|u_{0x}^2+\left(\kappa(\theta_0^\delta)\theta_{0x}\right)_x\right|}{Q^\prime(\theta_0^\delta)\sqrt{\rho_0^\delta}}
+|\sqrt{\rho_0^\delta}u_0\theta_{0x}|+|\sqrt{\rho_0^\delta}\theta_0^\delta
u_{0x}|\\&\le&\frac{\left|u_{0x}^2+\left(\kappa(\theta_0)\theta_{0x}\right)_x\right|}{Q^\prime(\theta_0^\delta)\sqrt{\rho_0^\delta}}
+\frac{\left|\left(\kappa(\theta_0^\delta)\theta_{0x}\right)_x-\left(\kappa(\theta_0)\theta_{0x}\right)_x\right|}{Q^\prime(\theta_0^\delta)\sqrt{\rho_0^\delta}}
+c\\&\le&c|g_2|+\frac{c\delta}{\sqrt{\rho_0^\delta}}(1+|\theta_{0xx}|)+c,
\eeqno which implies \beq\label{non-rho theta_t^2(0)}
\nonumber\int_I\rho\theta_t^2(0) & \le &
c\int_Ig_2^2+c\int_I\theta_{0xx}^2+c \\ & \le & c. \eeq Substituting
(\ref{non-rho theta_t^2(0)}) into (\ref{non-3.37}), using the
Gronwall inequality and Lemma \ref{non-le:3.8}, we complete the
proof.\qed
\begin{cor}\label{non-cor:3.2}Under the conditions of Theorem \ref{non-th:3.1}, it holds
$$\int_0^T\|\theta_t\|_{L^\infty}^2\le c.
$$
\end{cor}
\noindent{\it Proof.} By Lemma \ref{non-le:3.2}, ($A_5$), Corollary
\ref{non-cor:2.1}, Corollary \ref{non-cor:3.1}, Lemma
\ref{non-rle:3.11}, and $(\kappa\theta_t)_x=(\kappa\theta_x)_t$, we
get \beqno \int_0^T\|\kappa\theta_t\|_{L^\infty}^2 \le
c\int_0^T\|(\kappa\theta_t)_x\|_{L^2}^2 +c \le c. \eeqno This
combining ($A_5$) completes the proof.\qed
\begin{cor}\label{non-cor:3.3}Under the conditions of Theorem \ref{non-th:3.1}, it holds
$$\int_{Q_T}\theta_{xt}^2\le c.
$$
\end{cor}
\noindent{\it Proof.} Since
$$
\kappa\theta_{xt}=(\kappa\theta_x)_t-\kappa^\prime\theta_t\theta_x,
$$
we obtain \beqno
\int_{Q_T}\theta_{xt}^2&\le&c\int_{Q_T}\kappa^2\theta_{xt}^2\\&\le&c\int_{Q_T}\left|(\kappa\theta_x)_t\right|^2+c\int_{Q_T}(\kappa^\prime)^2\theta_t^2\theta_x^2
\\&\le&c+c\int_0^T\sup\limits_{x\in
I}\theta_t^2\int_I\theta_x^2\\&\le&c, \eeqno where we have used
($A_5$), Lemma \ref{non-le:3.8}, Lemma \ref{non-rle:3.11}, Corollary
\ref{non-cor:3.1} and Corollary \ref{non-cor:3.2}. \qed

\begin{cor}\label{non-rcor:3.4}Under the conditions of Theorem \ref{non-th:3.1}, it holds for any $0\le t\le T$
$$\|\theta\|_{W^{1,\infty}(Q_T)}+\int_I\theta_{xx}^2+\int_{Q_T}\theta_{xxx}^2\le
c.$$
\end{cor}
 \noindent{\it Proof.}  From (\ref{non-3.26}) and Lemma \ref{non-rle:3.11}, we have
 \beq\label{non-3.38}
\int_I\theta_{xx}^2\le c,
 \eeq which, combining Corollary \ref{non-cor:3.1}, Lemma
 \ref{non-le:3.8} and the Sobolev inequality, gives
 \be\label{non-3.41theta}
\|\theta\|_{W^{1,\infty}(Q_T)}\le c.
 \ee
Differentiating (\ref{non-3.6}) w.r.t. $x$, we have
\beq\label{non-3.39}\nonumber
&&\kappa\theta_{xxx}=-3\kappa^\prime\theta_x\theta_{xx}-\kappa^{\prime\prime}\theta_x^3-2u_xu_{xx}+\rho
Q^\prime\theta_{xt}+\rho_xQ^\prime\theta_t+\rho
Q^{\prime\prime}\theta_x\theta_t\\&&\ \ \ \ \ \ \ \ \ \ \ \ +(\rho u
Q^\prime\theta_x)_x+(\rho\theta Q^\prime u_x)_x. \eeq By
(\ref{non-3.38}), (\ref{non-3.41theta}), (\ref{non-3.39}), ($A_4$),
($A_5$), Lemma \ref{non-le:3.10}, Lemma \ref{non-rle:3.11}
 and Corollary \ref{non-rcor:3.3}, we have
\beq\label{non-3.40}\nonumber
\int_I\theta_{xxx}^2&\le&c\int_I\theta_x^2\theta_{xx}^2+c\int_I\theta_x^6+c\int_Iu_x^2u_{xx}^2+c\int_I\rho^2\theta_{xt}^2+c\int_I\rho_x^2\theta_t^2
+c\int_I\rho^2\theta_x^2\theta_t^2\\&&\nonumber+c\int_I\left|(\rho u
Q^\prime\theta_x)_x\right|^2+c\int_I\left|(\rho\theta Q^\prime
u_x)_x\right|^2+c\\&\le&\nonumber
c\int_I\rho^2\theta_{xt}^2+c\int_I\rho_x^2\theta_t^2
+c\\&\le&c\int_I\theta_{xt}^2+c\sup\limits_{x\in I}\theta_t^2+c.
\eeq By (\ref{non-3.40}), Corollary \ref{non-cor:3.2} and Corollary
\ref{non-cor:3.3}, we obtain
$$
\int_{Q_T}\theta_{xxx}^2\le c.
$$
\qed\\[0.5mm]

The next lemma, which we used in \cite{Ding-Wen-Zhu} to get $H^4-$estimates of velocity, plays an important role in getting $H^3-$estimates of $\theta$ in the following.
\begin{Lemma}\label{non-le:3.14} Under the conditions of Theorem \ref{non-th:3.1}, it holds
$$
\|(\sqrt{\rho})_x\|_{L^\infty(Q_T)}+\|(\sqrt{\rho})_t\|_{L^\infty(Q_T)}\le
c.
$$
\end{Lemma}
\noindent{\it Proof.} Multiplying $(\ref{non-1.2})_1$ by
$\D\frac{1}{2\sqrt{\rho}}$, we have \be\label{non-3.46}
(\sqrt{\rho})_t+(\sqrt{\rho})_xu+\frac{1}{2}\sqrt{\rho}u_x=0. \ee
Differentiating (\ref{non-3.46}) with respect to $x$, we get
$$
(\sqrt{\rho})_{xt}+(\sqrt{\rho})_{xx}u+\frac{3}{2}(\sqrt{\rho})_xu_x+
\frac{1}{2}\sqrt{\rho}u_{xx}=0.
$$
Denote $h=(\sqrt{\rho})_x$, we have
$$
h_t+h_x u+\frac{3}{2}hu_x+ \frac{1}{2}\sqrt{\rho}u_{xx}=0,
$$
which implies \be\label{non-3.47}
\frac{d}{ds}\left\{h\exp\left(\frac{3}{2}\int_0^s
\partial_Xu\left(X(\tau; x,
t),\tau\right)d\tau\right)\right\}=-\frac{1}{2}
\sqrt{\rho}(\partial_X^2u)\exp\left(\frac{3}{2}\int_0^s
\partial_Xu\left(X(\tau; x, t),\tau\right)d\tau\right), \ee where $X(s; x,
t)$ is the solution to (\ref{non-3.3}).

Integrating (\ref{non-3.47}) over $(0,t)$, we get
$$
\arraycolsep=1.5pt
\begin{array}{rl}
h(x,t)= & \displaystyle \exp\left(-\frac{3}{2}\int_0^t
\partial_Xu \left(X(\tau; x,
t),\tau\right)d\tau\right)h\left(X(0; x, t),0\right)
\\ [5mm]
& \displaystyle -\frac{1}{2}\exp\left(-\frac{3}{2}\int_0^t
\partial_Xu \left(X(\tau; x,
t),\tau\right)d\tau\right)
\int_0^t\sqrt{\rho}(\partial_X^2u)\exp\left(\frac{3}{2}\int_0^s
\partial_Xu\left(X(\tau; x, t),\tau\right)d\tau\right)ds.
\end{array}
$$
This together with Corollary \ref{non-rcor:3.3}, Lemma
\ref{non-le:3.10} and the Sobolev inequality, implies
\be\label{non-3.48} \|(\sqrt{\rho})_x\|_{L^\infty(Q_T)}\le c. \ee
From
 (\ref{non-3.46}), (\ref{non-3.48}), Lemma \ref{non-le:3.10} and Corollary \ref{non-rcor:3.3}, we get
$$
\|(\sqrt{\rho})_t\|_{L^\infty(Q_T)}\le c.
$$
This proves Lemma 3.10. \qed\\[1.5mm]

 The next lemma will be used to get $H^3-$estimates of $\theta$.
\begin{Lemma}\label{non-le:3.15}Under the conditions of Theorem \ref{non-th:3.1}, it holds for any $0\le t\le T$
$$
\int_I\rho^2\left|(\kappa\theta_x)_t\right|^2+\int_{Q_T}\rho^3\theta_{tt}^2\le
c.
$$
\end{Lemma}
\noindent{\it Proof.} Multiply $(\ref{non-3.36})$ by
$\rho^{\gamma_1}(\kappa\theta_t)_t$ (i.e.
$\rho^{\gamma_1}\kappa\theta_{tt}+\rho^{\gamma_1}\kappa^\prime\theta_t^2$,
where $\gamma_1$ is to be decided later), and using integration by
parts, we have
 \beq\label{non-r3.41}
\nonumber&&\int_I\rho^{\gamma_1+1}\kappa
Q^\prime\theta_{tt}^2+\frac{1}{2}\frac{d}{dt}\int_I\rho^{\gamma_1}\left|(\kappa\theta_x)_t\right|^2\\&=&\nonumber\frac{\gamma_1}{2}\int_I\rho^{\gamma_1-1}\rho_t\left|(\kappa\theta_x)_t\right|^2
-\gamma_1\int_I\rho^{\gamma_1-1}\rho_x\kappa\theta_{tt}(\kappa\theta_x)_t-\gamma_1\int_I\rho^{\gamma_1-1}\rho_x\kappa^\prime\theta_t^2(\kappa\theta_x)_t
\\&&\nonumber+2\int_Iu_xu_{xt}(\rho^{\gamma_1}\kappa\theta_{tt}+\rho^{\gamma_1}\kappa^\prime\theta_t^2)-\int_I\rho^{\gamma_1+1}Q^\prime\kappa^\prime\theta_t^2\theta_{tt}\\&&
-\int_I\left(\rho
Q^{\prime\prime}\theta_t^2+\rho_tQ^\prime\theta_t+(\rho
uQ^\prime\theta_x)_t+(\rho\theta Q^\prime
u_x)_t\right)\left(\rho^{\gamma_1}\kappa\theta_{tt}+\rho^{\gamma_1}\kappa^\prime\theta_t^2\right).
\eeq We are going to look for the minimal of $\gamma_1$. It seems
that the second term of the right side plays an important role.
\beq\label{non-r3.42}
-\gamma_1\int_I\rho^{\gamma_1-1}\rho_x\kappa\theta_{tt}(\kappa\theta_x)_t&=&\nonumber
-2\gamma_1\int_I\rho^{\gamma_1-\frac{1}{2}}(\sqrt{\rho})_x\kappa\theta_{tt}(\kappa\theta_x)_t\\&\le&\frac{1}{4}\int_I\rho^{\gamma_1+1}\kappa
Q^\prime\theta_{tt}^2+c\int_I\rho^{\gamma_1-2}|(\kappa\theta_x)_t|^2,
\eeq where we have used Lemma \ref{non-le:3.14}, ($A_4$), ($A_5$),
Corollary \ref{non-cor:3.1} and the Cauchy inequality.

From Lemma \ref{non-rle:3.11}, we know that
$\int_{Q_T}|(\kappa\theta_x)_t|^2\le c$. This implies that the
minimal of $\gamma_1$ should be 2. Substituting $\gamma_1=2$ into
(\ref{non-r3.41}) and (\ref{non-r3.42}), and then substituting
(\ref{non-r3.42}) into (\ref{non-r3.41}), we have \beqno
\nonumber&&\frac{3}{4}\int_I\rho^3\kappa
Q^\prime\theta_{tt}^2+\frac{1}{2}\frac{d}{dt}\int_I\rho^2\left|(\kappa\theta_x)_t\right|^2\\&\le&
c\int_I\left|(\kappa\theta_x)_t\right|^2+c\int_I\rho\theta_t^4
+c\int_Iu_{xt}^2+c\int_I\theta_{xt}^2+c\\&&+c\left(\int_I
\rho^3\kappa
Q^\prime\theta_{tt}^2\right)^\frac{1}{2}\left\{1+\left(\int_I
\rho\theta_t^4\right)^\frac{1}{2}+\|\theta_{xt}\|_{L^2}+\|\sqrt{\rho}u_t\|_{L^2}+\|u_{xt}\|_{L^2}+\|\sqrt{\rho}\theta_t\|_{L^2}\right\}
\\&\le&\frac{1}{4}\int_I\rho^3\kappa
Q^\prime\theta_{tt}^2+c\int_I\left|(\kappa\theta_x)_t\right|^2+c\|\theta_t\|_{L^\infty}^2+c\int_I\theta_{xt}^2+c\int_Iu_{xt}^2+c,\eeqno
where we have used ($A_4$), ($A_5$), Lemma \ref{non-le:3.8}, Lemma
\ref{non-le:3.10}, Lemma \ref{non-rle:3.11}, Corollary
\ref{non-rcor:3.3}, Corollary \ref{non-rcor:3.4} and the Cauchy
inequality. This implies \beqno \int_I\rho^3\kappa
Q^\prime\theta_{tt}^2+\frac{d}{dt}\int_I\rho^2\left|(\kappa\theta_x)_t\right|^2\le
c\int_I\left|(\kappa\theta_x)_t\right|^2+c\|\theta_t\|_{L^\infty}^2+c\int_I\theta_{xt}^2+c\int_Iu_{xt}^2+c.
\eeqno Integrating it over $(0, t)$, and using ($A_4$), ($A_5$),
Lemma \ref{non-le:3.8}, Lemma \ref{non-rle:3.11}, Corollary
\ref{non-cor:3.2} and Corollary \ref{non-cor:3.3}, we have
\be\label{non-3.54}
\int_I\rho^2\left|(\kappa\theta_x)_t\right|^2+\int_0^t\int_I\rho^3\theta_{tt}^2\le
c\int_I\rho^2\left|(\kappa\theta_x)_t\right|^2(0)+c. \ee By ($A_5$),
(\ref{non-rho theta_t^2(0)}), (\ref{non-3.39}) and Corollary
\ref{non-cor:3.1}, we get \beq\label{non-rhotheta{xt}}\nonumber
\int_I\rho^2\left|(\kappa\theta_x)_t\right|^2(0)&\le&c\int_I\rho^2\theta_{xt}^2(0)+c\int_I\rho\theta_t^2(0)\\&\le&\nonumber
c\|\theta_0\|_{H^3}^2
+c\|u_0\|_{H^2}^2+c\int_I\rho_x^2u_t^2(0)+c\\&\le&\nonumber
c+c\int_I|(\sqrt{\rho})_x|^2\rho u_t^2(0)\\&\le&c. \eeq Substituting
(\ref{non-rhotheta{xt}}) into (\ref{non-3.54}), we complete the
proof. \qed
\begin{cor}\label{non-cor:3.4}Under the conditions of Theorem \ref{non-th:3.1}, it holds for any $0\le t\le T$
$$
\int_I\left(\theta_{xxx}^2+\rho^2\theta_{xt}^2\right)\le c.
$$
\end{cor}
\noindent{\it Proof.} A direct calculation gives
$$
\rho\kappa\theta_{xt}=\rho(\kappa\theta_x)_t-\rho\kappa^\prime\theta_t\theta_x,
$$
which implies \beq\label{non-r3.51}
\nonumber\int_I\rho^2\theta_{xt}^2&\le&c\int_I\rho^2\left|(\kappa\theta_x)_t\right|^2+c\|\theta_x\|_{L^\infty}^2\int_I\rho\theta_t^2 \\
& \le & c. \eeq
Here we have used $(A_5)$, Lemma \ref{non-le:3.10},
Lemma \ref{non-rle:3.11}, Lemma \ref{non-le:3.15} and Corollary
\ref{non-rcor:3.4}.

 From the second inequality of (\ref{non-3.40}),
we obtain \beqno\int_I\theta_{xxx}^2&\le&
c\int_I\rho^2\theta_{xt}^2+c\int_I\rho_x^2\theta_t^2
+c\\&\le&c\int_I|(\sqrt{\rho})_x|^2\rho\theta_t^2 +c\\&\le&c,\eeqno
where we have used Lemma \ref{non-rle:3.11}, (\ref{non-r3.51}) and
Lemma \ref{non-le:3.14}.\qed\\[1.5mm]

 The next lemma will be used to get $H^3-$estimates of $u$.

\begin{Lemma}\label{non-rle:3.12}Under the conditions of Theorem \ref{non-th:3.1}, it holds for any $0\le t\le T$
$$
\int_I\rho^2u_{xt}^2+\int_{Q_T}\rho^3u_{tt}^2\le c.
$$
\end{Lemma}
\noindent{\it Proof.} Similarly to Lemma \ref{non-le:3.15},
multiplying (\ref{non-3.22}) by $\rho^2u_{tt}$, and integrating it
over $I$, we have \beqno
&&\int_I\rho^3u_{tt}^2+\frac{1}{2}\frac{d}{dt}\int_I\rho^2
u_{xt}^2\\&=&\int_I\rho\rho_tu_{xt}^2-2\int_I\rho\rho_xu_{xt}u_{tt}-\int_I\rho^2u_{tt}(\rho_tu_t+\rho_tuu_x+\rho
u_tu_x+\rho
uu_{xt}+P_{xt})\\&\le&c\int_Iu_{xt}^2-4\int_I\rho^{\frac{3}{2}}(\sqrt{\rho})_xu_{xt}u_{tt}+\frac{1}{4}\int_I\rho^3u_{tt}^2+c\int_I\rho
u_t^2+c\int_I|(\rho
Q)_{xt}|^2+c\int_I|(P_c)_{xt}|^2+c\\&\le&\frac{1}{2}\int_I\rho^3u_{tt}^2+c\int_Iu_{xt}^2+c\|\theta_t\|_{L^\infty}^2+c\int_I\theta_{xt}^2+c.
\eeqno Here, we have used integration by parts, Lemma
\ref{non-le:3.8}, Lemma \ref{non-le:3.10}, Lemma \ref{non-le:3.14},
Corollary \ref{non-rcor:3.3}, Corollary \ref{non-rcor:3.4} and the
Cauchy inequality.

The first term of the right side can be absorbed by the left. After
that, we have \beqno \int_I\rho^3u_{tt}^2+\frac{d}{dt}\int_I\rho^2
u_{xt}^2\le
c\int_Iu_{xt}^2+c\|\theta_t\|_{L^\infty}^2+c\int_I\theta_{xt}^2+c.
\eeqno Integrating this inequality on both side over $(0,t)$, and
using Lemma \ref{non-le:3.8}, Corollary \ref{non-cor:3.2} and
Corollary \ref{non-cor:3.3}, we have \be\label{non-r3.47}
\int_0^t\int_I\rho^3u_{tt}^2+\int_I\rho^2 u_{xt}^2\le \int_I\rho^2
u_{xt}^2(0)+c. \ee Similarly to (\ref{non-rhotheta{xt}}), we use
(\ref{non-rhou_t^2(0)}), (\ref{non-3.29}), ($A_3$) and ($A_4$) to
get \beq\label{non-r3.49} \nonumber\int_I\rho^2u_{xt}^2(0)&\le&
c\|u_0\|_{H^3}^2+c\|\theta_0\|_{H^2}^2+c\|\rho_0\|_{H^2}^2+c\int_I\rho
u_t^2(0)+c\\&\le&c. \eeq
Substituting (\ref{non-r3.49}) into (\ref{non-r3.47}), we complete the proof.\qed\\[0.8mm]

By (\ref{non-3.30}), Lemma \ref{non-le:3.8}, Lemma
\ref{non-le:3.10}, Lemma \ref{non-le:3.14}, Lemma \ref{non-rle:3.12}
and Corollary \ref{non-rcor:3.4}, we get the following corollary.
\begin{cor}\label{non-rcor:3.8}Under the conditions of Theorem \ref{non-th:3.1}, it holds for any $0\le t\le T$
$$\int_Iu_{xxx}^3\le c.$$
\end{cor}

\bigskip

From the above estimates, we get \beq\label{non-r3.54}\nonumber
&&\|(\sqrt{\rho})_x\|_{L^\infty}+\|(\sqrt{\rho})_t\|_{L^\infty}+\|\rho\|_{H^2}
+\|\rho_t\|_{H^1}+
 \displaystyle \|u\|_{H^3}
+\|\rho
u_t\|_{H^1}+\|\sqrt{\rho}u_t\|_{L^2}+\|\theta\|_{H^3}\\&&+\|\sqrt{\rho}\theta_t\|_{L^2}+\|\rho\theta_{t}\|_{H^1}
+\int_{Q_T}\left( u_{xt}^2+\rho_{tt}^2
+\theta_t^2+\theta_{xt}^2+\rho^3u_{tt}^2+\rho^3\theta_{tt}^2\right)\le
c.\eeq

\begin{cor}\label{non-cor:3.9} Under the conditions of Theorem \ref{non-th:3.1}, there exists a positive constant $c_\delta$ depending on $\delta$ such that for any $(x,t)\in Q_T$,  it holds
\be\label{non-r3.55}
\begin{cases}\rho(x,t)\ge\D\frac{\delta}{c}>0,\\
\theta(x,t)\ge c_\delta>0.
\end{cases}\ee
\end{cor}
\noindent{\it Proof.} By (\ref{non-r3.5}), ($A_3$), ($A_4$), Lemma
\ref{non-le:3.10} and Corollary \ref{non-rcor:3.4}, we have for any
$(x,t)\in Q_T$ \beqno \rho(x,t)\ge\frac{\delta}{c}. \eeqno This gets
(\ref{non-r3.55})$_1$. (\ref{non-r3.55})$_2$ can be got by
(\ref{non-r3.55})$_1$, (\ref{non-r3.54}), (\ref{non-3.6}) and
the maximum principle for parabolic equation.\qed\\[1.5mm]

From (\ref{non-r3.54}), (\ref{non-r3.55}), (\ref{non-3.22}) and
(\ref{non-3.36}), we obtain \beqno &&\|\rho\|_{H^2}
+\|\rho_t\|_{H^1}+
 \|u\|_{H^3}
+\|u_t\|_{H^1}+\|\theta\|_{H^3}+\|\theta_{t}\|_{H^1} \\&&
+\int_{Q_T}\left( u_{xt}^2+u_{xxt}^2+\rho_{tt}^2
+\theta_t^2+\theta_{xt}^2+\theta_{xxt}^2+u_{tt}^2+\theta_{tt}^2\right)\le
c. \eeqno This proves Theorem \ref{non-th:3.1}.\qed\\[1.5mm]

{\noindent\bf Proof of Theorem \ref{non-rth:1.1}:}

Consider (\ref{non-1.2})-(\ref{non-1.4}) with initial data replaced
by ($\rho_0^\delta$, $u_0$, $\theta_0^\delta$), we obtain from
Theorem \ref{non-th:3.1} that there exists a unique solution
($\rho^{\delta}$, $u^{\delta}$, $\theta^{\delta}$), such that
(\ref{non-r3.54}) and (\ref{non-r3.55}) are valid when we replace
($\rho,$ $u$, $\theta$) by ($\rho^{\delta}$, $u^{\delta}$,
$\theta^{\delta}$). With the estimates uniform for $\delta$, we take
$\delta\rightarrow0^+$ (take subsequence if necessary) to get a
solution to (\ref{non-1.2})-(\ref{non-1.4}) still denoted by
($\rho$, $u$, $\theta$) which satisfies (\ref{non-r3.54}) by the
lower semi-continuity of the norms. This proves the existence of the
solutions as in Theorem \ref{non-rth:1.1}. The uniqueness of the
solutions can be proved by the standard method like in
\cite{cho-Kim: perfect gas}, we omit it for brevity. The proof of
Theorem \ref{non-rth:1.1} is complete.\qed

 \setcounter{section}{3}
\setcounter{equation}{0}
\section{ \ Proof of Theorem \ref{non-th:1.1}}
In this section, we use the similar arguments as in Section 3 to
prove Theorem \ref{non-th:1.1}.
 Throughout this section, we denote $c$ to
be a generic constant depending on $\rho_0$, $u_0$, $\theta_0$, $T$
and some other known constants but independent of $\delta$ for any
$\delta\in(0,1)$.

Denote $\rho_0^\delta=\rho_0+\delta$,
$\theta_0^\delta=\theta_0+\delta$ and $P^\delta_0=P(\rho_0^\delta,\
\theta_0^\delta)$, where $\rho_0$ and $\theta_0$ satisfy the same
conditions as those in Theorem \ref{non-th:1.1}. Note that
$\rho_0^\delta\in H^4(I)$, $\rho_0^\delta\ge\delta>0$,
$\theta_0^\delta\in H^3(I)$,
$\partial_x\theta_0^\delta|_{x=0,1}=\partial_x\theta_0|_{x=0,1}=0$,
and \be\label{non-4.1}
\begin{cases}\|\rho_0^\delta\|_{H^4}\le c,\\ \|\big(\sqrt{\rho_0^\delta}\big)_x\|_{L^\infty}\le
c,\\ \|\theta_0^\delta\|_{H^3}\le c.
\end{cases}\ee
Different from Section 3, we need to mollify $g_3$. Denote
$g_3^\delta=J_\delta*\overline{g}_3$, then $g_3^\delta\in
C^\infty(I)$, where
$$\overline{g}_3(x)=\begin{cases} -g_3(-x),\ x\in [-1, 0),\\ g_3(x),\ \ \ \ \ x\in I,\\
-g_3(2-x),\ x\in (1,2],
\end{cases}$$
and
$J_\delta(\cdot)=\frac{1}{\sqrt{\delta}}J(\frac{\cdot}{\sqrt{\delta}})$,
and $J$ is the usual mollifier such that $J\in
C_0^\infty(\mathbb{R})$, supp$J\in(-1, 1)$, and
$\int_\mathbb{R}J(x)dx=1$. Since $g_3\in H_0^1(I)$, we have
$\overline{g}_3\in H_0^1([-1,2])$ and
$$\partial_x\overline{g}_3(x)=\begin{cases}g_3^\prime(-x),\ x\in [-1, 0),\\ g_3^\prime(x),\ \ \ \  x\in I,\\
g_3^\prime(2-x),\ x\in (1, 2].
\end{cases}$$
{\it Claim:} \beq\label{non-r3.1}\begin{cases} g_3^\delta\rightarrow
g_3\ \ \mathrm{in}\ H^1(I),\ \mathrm{as}\
\delta\rightarrow0,\\
\|g_3^\delta\|_{H^1(I)}\le c\|\overline{g}_3\|_{H^1([-1,2])}\le
c\|g_3\|_{H^1(I)},\ \mathrm{for}\ \mathrm{any}\ \delta\in(0, 1),
\\
\|\sqrt{\rho^\delta_0}(g_3^\delta)_{xx}\|_{L^2(I)}\le c,\
\mathrm{for}\ \mathrm{any}\ \delta\in(0, 1).\end{cases} \eeq In
fact, the proof of (\ref{non-r3.1})$_1$ and (\ref{non-r3.1})$_2$ can
be found in \cite{26}. We are going to prove (\ref{non-r3.1})$_3$.
\beq\label{non-4.3}
\nonumber\sqrt{\rho_0^\delta}(g_3^\delta)_{xx}&=&(\sqrt{\rho_0^\delta}-\sqrt{\rho_0})(g_3^\delta)_{xx}+\sqrt{\rho_0}(g_3^\delta)_{xx}\\&=&
\nonumber\frac{\delta(g_3^\delta)_{xx}}{\sqrt{\rho_0^\delta}+\sqrt{\rho}}+\sqrt{\rho_0}(g_3^\delta)_{xx}\\&=&A_1+A_2.
\eeq Recall
$J_\delta(\cdot)=\frac{1}{\sqrt{\delta}}J(\frac{\cdot}{\sqrt{\delta}})$,
we conclude \beq\label{non-4.4}\nonumber
\|A_1\|_{L^2(I)}&\le&\sqrt{\delta}\|(g_3^\delta)_{xx}\|_{L^2(I)}\\&\le&\nonumber
c\|(\overline{g}_3)_x\|_{L^2([-1,2])}\\&\le&c\|(g_3)_x\|_{L^2(I)}.
\eeq A direct calculation combining
$\left(\sqrt{\rho_0}(g_3)_x\right)_x\in L^2(I)$ gives
\beq\label{non-4.5} \int_I|A_2|^2\le c. \eeq By (\ref{non-4.3}),
(\ref{non-4.4}) and (\ref{non-4.5}), we get (\ref{non-r3.1})$_3$.

 Let
$u_0^\delta$ be the solution to the following elliptic problem for
each $\delta\in(0,1)$: \beq\label{non-r3.2}\begin{cases} u^\delta_{0xx}-(P_0^\delta)_x=\rho^\delta_0g^\delta_3,\\
u_0^\delta|_{x=0,1}=0.\end{cases}
 \eeq
Since $\rho_0^\delta=\rho_0+\delta\in H^4(I)$,
$\theta_0^\delta=\theta_0+\delta\in H^3(I)$, and $g_3^\delta\in
C^\infty(I)$, we obtain from the elliptic theory (see \cite{26}),
(\ref{non-4.1}), (\ref{non-r3.1}) and (\ref{non-r3.2}) that
$u_0^\delta\in H^4(I)\cap H^1_0(I)$ with the following properties:
\beq\label{non-4.6}
\begin{cases}
u_0^\delta\rightarrow u_0\ \mathrm{in}\ H^3(I),\
\mathrm{as}\ \delta\rightarrow0,\\
\|u_0^\delta\|_{H^4(I)}\le c\ \ \mathrm{for}\ \mathrm{any}\
\delta\in(0, 1).
\end{cases}
\eeq

\begin{Theorem}\label{non-th:4.1} Consider the same assumptions as in Theorem \ref{non-th:1.1}. Then for any $T>0$ and $\delta\in(0,1)$ there exists a unique global solution $(\rho,
u, \theta)$ to (\ref{non-1.2})-(\ref{non-1.4}) with initial data
replaced by ($\rho_0^\delta,u_0^\delta,\theta_0^\delta$), such that
 \beqno &\rho\in C([0,T];H^4), \ \ \ \rho_t\in
C([0,T];H^3), \ \ \ \rho_{tt}\in C([0,T];H^1)\cap L^2([0,T];H^2),
&\\&\rho_{ttt}\in L^2(Q_T),\ \rho\ge\frac{\D\delta}{c}>0,\ u\in
C([0,T];H^4\cap H^1_0)\cap L^2([0,T];H^5),&\\& u_t\in
C([0,T];H^2)\cap L^2([0,T];H^3), \ \ \ u_{tt}\in C([0,T];L^2)\cap
L^2([0,T];H_0^1),\ &\\& \theta\in C([0,T]; H^3)\cap L^2([0,T]; H^4),
\ \theta_t\in C([0,T]; H^1)\cap L^2([0,T]; H^2),\ &\\&
\theta_{tt}\in L^2([0,T]; L^2),\ \theta\ge c_\delta>0,\ & \eeqno
where $c_\delta$ is a constant depending on $\delta$, but
independent of $u$.
\end{Theorem}

{\noindent\bf Proof of Theorem \ref{non-th:4.1}:}

Similarly to the proof of Theorem \ref{non-th:3.1}, Theorem
\ref{non-th:4.1} can be proved by some {\it a priori} estimates
globally in time.

For any given $T\in(0,+\infty)$, let $(\rho,u,\theta)$ be the
solution to (\ref{non-1.2})-(\ref{non-1.4}) as in Theorem
\ref{non-th:4.1}. Then we have the following estimates.

\begin{Lemma}\label{non-le:4.1} Under the conditions of Theorem \ref{non-th:4.1}, it holds for any $0\le t\le T$
\beqno
&&\|(\sqrt{\rho})_x\|_{L^\infty}+\|(\sqrt{\rho})_t\|_{L^\infty}+\|\rho\|_{H^2}
+\|\rho_t\|_{H^1}+
 \displaystyle \|u\|_{H^3}
+\|\rho
u_t\|_{H^1}+\|\sqrt{\rho}u_t\|_{L^2}+\|\theta\|_{H^3}\\&&+\|\sqrt{\rho}\theta_t\|_{L^2}+\|\rho\theta_{t}\|_{H^1}
+\int_{Q_T}\left( u_{xt}^2+\rho_{tt}^2
+\theta_t^2+\theta_{xt}^2+\rho^3u_{tt}^2+\rho^3\theta_{tt}^2\right)\le
c.\eeqno
\end{Lemma}
\noindent{\it Proof.} Though the initial velocity in Theorem
\ref{non-th:4.1} (i.e. $u_0^\delta$) is different from that in
Theorem \ref{non-th:3.1} (i.e. $u_0$), both of them are bounded in
$H^3$. It suffices to check if (\ref{non-rhou_t^2(0)}) and
(\ref{non-rho theta_t^2(0)}) work here.  If do, Lemma \ref{non-le:4.1} will be obtained from (\ref{non-r3.54}).

By (\ref{non-3.10}) and
(\ref{non-r3.2}) \beqno |\sqrt{\rho}
u_t(x,0)|&\le&\frac{\left|u_{0xx}^\delta-P(\rho_0^\delta,\theta_0^\delta)_x\right|}{\sqrt{\rho_0^\delta}}+\sqrt{\rho_0^\delta}
|u_0^\delta
u_{0x}^\delta|\\&=&\sqrt{\rho_0^\delta}|g_3^\delta|+\sqrt{\rho_0^\delta}
|u_0^\delta u_{0x}^\delta|. \eeqno This gives
$$
\int_I\rho u_t^2(0)\le c.
$$
Therefore, (\ref{non-rhou_t^2(0)}) is valid here.

Multiplying (\ref{non-3.6}) by $\D\frac{1}{
Q^\prime(\theta)\sqrt{\rho}}$, taking $t\rightarrow0^+$, and using
(\ref{non-1.5})$_2$, we have \beqno
|\sqrt{\rho}\theta_t(x,0)|&\le&\frac{\left|(u^\delta_{0x})^2+\left(\kappa(\theta_0^\delta)\theta_{0x}\right)_x\right|}{Q^\prime(\theta_0^\delta)\sqrt{\rho_0^\delta}}
+|\sqrt{\rho_0^\delta}u_0^\delta\theta_{0x}|+|\sqrt{\rho_0^\delta}\theta_0^\delta
u_{0x}^\delta|\\&\le&\frac{\left|u_{0x}^2+\left(\kappa(\theta_0)\theta_{0x}\right)_x\right|}{Q^\prime(\theta_0^\delta)\sqrt{\rho_0^\delta}}
+\frac{c\left|u_{0x}^\delta-u_{0x}\right|}{Q^\prime(\theta_0^\delta)\sqrt{\rho_0^\delta}}
+\frac{\left|\left(\kappa(\theta_0^\delta)\theta_{0x}\right)_x-\left(\kappa(\theta_0)\theta_{0x}\right)_x\right|}{Q^\prime(\theta_0^\delta)\sqrt{\rho_0^\delta}}
+c\\&\le&c|g_2|+\frac{c\delta}{\sqrt{\rho_0^\delta}}(1+|\theta_{0xx}|)+\frac{c\left|u_{0x}^\delta-u_{0x}\right|}{\sqrt{\delta}}.
\eeqno Note that $\|u_{0x}^\delta-u_{0x}\|_{L^2(I)}\le
c\sqrt{\delta}$ by (\ref{non-1.5})$_1$ and (\ref{non-r3.2}). This
gives \beqno\int_I\rho\theta_t^2(0) \le
c\int_Ig_2^2+c\int_I\theta_{0xx}^2+c \le c. \eeqno Therefore,
(\ref{non-rho theta_t^2(0)}) is valid here. \qed
\begin{Lemma}\label{non-le:4.2} Under the conditions of Theorem
\ref{non-th:4.1}, it holds for any $0\le t\le T$
$$
\int_Iu_{xt}^2+\int_{Q_T}\rho u_{tt}^2\le c.
$$
\end{Lemma}
\noindent{\it Proof.} Multiplying (\ref{non-3.22}) by $u_{tt}$,
integrating it over $I$, and using integration by parts, Lemma
\ref{non-le:2.2}, Lemma \ref{non-le:4.1} and the Cauchy inequality,
we have \beqno &&\int_I\rho
u_{tt}^2+\frac{1}{2}\frac{d}{dt}\int_Iu_{xt}^2
\\&=&-\frac{1}{2}\frac{d}{dt}\int_I\rho_tu_t^2+\frac{1}{2}\int_I\rho_{tt}u_t^2-\frac{d}{dt}\int_I\rho_tuu_xu_t+\int_I\rho_{tt}uu_xu_t+\int_I\rho_tu_t^2u_x
\\&&+\int_I\rho_tuu_{xt}u_t-\int_I\rho u_tu_xu_{tt}-\int_I\rho
uu_{xt}u_{tt}+\frac{d}{dt}\int_IP_tu_{xt}-\int_IP_{tt}u_{xt}\\&\le&\frac{d}{dt}
\int_I\left(P_tu_{xt}-\frac{1}{2}\rho_tu_t^2-\rho_tuu_xu_t\right)+c\int_Iu_{xt}^2\int_I\rho_{tt}^2\\&&
+c\int_Iu_{xt}^2+c\int\rho_{tt}^2+\frac{1}2\int_I\rho
u_{tt}^2-\int_IP_{tt}u_{xt}+c. \eeqno This gives
\beq\label{non-3.32} \nonumber \int_I\rho
u_{tt}^2+\frac{d}{dt}\int_Iu_{xt}^2 &\le&\frac{d}{dt}
\int_I\left(2P_tu_{xt}-\rho_tu_t^2-2\rho_tuu_xu_t\right)+c\int_Iu_{xt}^2\int_I\rho_{tt}^2+c\int_Iu_{xt}^2\\&&+c\int\rho_{tt}^2
-\frac{d}{dt}\int_IP_t^2-2\int_IP_{tt}(u_{xt}-P_t)+c. \eeq We are
going to estimate the last term of the right side of
(\ref{non-3.32}). By ($A_2$)-($A_4$),  integration by parts, Lemma
\ref{non-le:4.1} and the Cauchy inequality, we have\beqno
&&-2\int_IP_{tt}(u_{xt}-P_t)\\&=&-2\int_I(\rho
Q)_{tt}(u_{xt}-P_t)-2\int_I(P_c)_{tt}\left[u_{xt}-(\rho
Q)_t-(P_c)_t\right] \\ & \le
&-2\int_I\left[(\kappa\theta_x)_x+u_x^2-(\rho u Q)_x-\rho\theta
Q^\prime u_x\right]_t(u_{xt}-P_t)
\\ && +c\int_I\rho_{tt}^2+c\int_Iu_{xt}^2+c\int_I\rho\theta_t^2+c
\\&=&2\int_I(\kappa\theta_x)_t(u_{xx}-P_x)_t-4\int_Iu_xu_{xt}(u_{xt}-P_t)
-2\int_I(\rho uQ)_t(u_{xx}-P_x)_t\\&&+2\int_I\left[\rho\theta Q^\prime
u_x\right]_t(u_{xt}-P_t)+c\int_I\rho_{tt}^2+c\int_Iu_{xt}^2+c\int_I\rho\theta_t^2+c.
\eeqno
This, combining (\ref{non-3.10}), ($A_2$)--($A_4$), Lemma \ref{non-le:4.1} and the
Cauchy inequality, concludes
\beq\label{non-3.33}\nonumber
&&-2\int_IP_{tt}(u_{xt}-P_t)
\\&\le & \nonumber c+2\int_I(\kappa\theta_x)_t(\rho
u_t+\rho uu_x)_t+c\int_Iu_{xt}^2+c\int_I\rho\theta_t^2 \\ &&
\nonumber -2\int_I(\rho uQ)_t(\rho u_t+\rho
uu_x)_t+2\int_I(\rho\theta Q^\prime
u_x)_t(u_{xt}-P_t)+c\int_I\rho_{tt}^2\\&\le&
c\int_I|(\kappa\theta_x)_t|^2+\frac{1}{2}\int_I\rho
u_{tt}^2+c\int_Iu_{xt}^2+c\int_I\rho\theta_t^2+c\int_I\rho_{tt}^2+c.
\eeq Substituting (\ref{non-3.33}) into (\ref{non-3.32}), we get
\beq\label{non-3.34} \nonumber \frac{1}{2}\int_I\rho
u_{tt}^2+\frac{d}{dt}\int_Iu_{xt}^2 &\le&\frac{d}{dt}
\int_I\left(2P_tu_{xt}-\rho_tu_t^2-2\rho_tuu_xu_t\right)+c\int_Iu_{xt}^2\int_I\rho_{tt}^2+c\int_Iu_{xt}^2\\&&+c\int\rho_{tt}^2
-\frac{d}{dt}\int_IP_t^2+c\int_I|(\kappa\theta_x)_t|^2+c\int_I\rho\theta_t^2+c.
\eeq Integrating (\ref{non-3.34}) over $(0, t)$, and using
(\ref{non-1.2})$_1$, integration by parts, (\ref{non-3.7}),
(\ref{non-3.10}), (\ref{non-r3.1})$_2$, (\ref{non-r3.2}), and Lemma
\ref{non-le:4.1}, we have
\beqno &&\frac{1}{2}\int_0^t\int_I\rho u_{tt}^2+\int_Iu_{xt}^2\\
&\le& \int_I\left(2P_tu_{xt}+(\rho
u)_xu_t^2-2\rho_tuu_xu_t\right)+c\int_0^t\int_Iu_{xt}^2\int_I\rho_{tt}^2
+c\int_0^t\int_I|(\kappa\theta_x)_t|^2+c\\&=&
\int_I\left(2P_tu_{xt}-2\rho u
u_tu_{xt}-2\rho_tuu_xu_t\right)+c\int_0^t\int_Iu_{xt}^2\int_I\rho_{tt}^2
+c\int_0^t\int_I|(\kappa\theta_x)_t|^2+c\\&\le&\frac{1}{2}\int_Iu_{xt}^2+c\int_I\rho\theta_t^2+c\int_I\rho^2u^2u_t^2+c\int_I\rho_t^2u^2u_x^2+c\int_0^t\int_Iu_{xt}^2\int_I\rho_{tt}^2
+c\int_0^t\int_I|(\kappa\theta_x)_t|^2+c\\&\le&\frac{1}{2}\int_Iu_{xt}^2+c\int_I\rho\theta_t^2++c\int_0^t\int_Iu_{xt}^2\int_I\rho_{tt}^2
+c\int_0^t\int_I|(\kappa\theta_x)_t|^2+c, \eeqno which implies
\beq\label{non-3.35} \int_0^t\int_I\rho u_{tt}^2+\int_Iu_{xt}^2 \le
c\int_I\rho\theta_t^2+c\int_0^t\int_Iu_{xt}^2\int_I\rho_{tt}^2
+c\int_0^t\int_I|(\kappa\theta_x)_t|^2+c.\eeq Using the Gronwall
inequality and Lemma \ref{non-le:4.1}, we complete the proof of the
lemma. \qed

\begin{cor}\label{non-cor:4.1}Under the conditions of Theorem \ref{non-th:4.1}, it holds for any $0\le t\le T$
$$\int_{Q_T}u_{xxt}^2\le
c.$$
\end{cor}
 \noindent{\it Proof.}
It follows from (\ref{non-3.22}), Lemma \ref{non-le:4.1} and
($A_2$)--($A_4$) that \beqno \int_{Q_T}u_{xxt}^2 & \le &
c\int_{Q_T}\rho
u_{tt}^2+c\int_{Q_T}\rho_t^2u_t^2+c\int_{Q_T}\rho_t^2u^2u_x^2+c\int_{Q_T}\rho^2u_t^2u_x^2
\\ && +c\int_{Q_T}\rho^2u^2u_{xt}^2+c\int_{Q_T}\left|(\rho
Q)_{xt}\right|^2+c\int_{Q_T}\left|(P_c)_{xt}\right|^2\\&\le&c+c\int_0^T\|u_t\|_{L^\infty}^2+c\int_{Q_T}u_{xt}^2+c\int_{Q_T}\rho_{xt}^2
+c\int_0^T\|\theta_t\|_{L^\infty}^2+c\int_{Q_T}\theta_{xt}^2+c\\&\le&c.
\eeqno This proves Corollary 4.1. \qed

\begin{Lemma}\label{non-le:4.3}Under the conditions of Theorem \ref{non-th:4.1}, it holds for any $0\le t\le T$
$$
\int_I\left(\rho_{xxx}^2+\rho_{xxt}^2+\rho_{tt}^2\right)+\int_{Q_T}\left(\rho_{xtt}^2+u_{xxxx}^2\right)\le
c.
$$
\end{Lemma}
\noindent{\it Proof.} Differentiating (\ref{non-3.27}) with respect
to $x$, we have \be\label{non-3.41}
\rho_{xxxt}=-\rho_{xxxx}u-4\rho_{xxx}u_x-6\rho_{xx}u_{xx}-4\rho_xu_{xxx}-\rho
u_{xxxx}. \ee Multiplying (\ref{non-3.41}) by $2\rho_{xxx}$,
integrating the resulting equation over $I$, and using integration
by parts and the H\"older inequality, we have \beqno
\frac{d}{dt}\int_I\rho_{xxx}^2&=&-7
\int_I\rho_{xxx}^2u_x-12\int_I\rho_{xx}\rho_{xxx}u_{xx}
-8\int_I\rho_x\rho_{xxx}u_{xxx}-2\int_I\rho\rho_{xxx}u_{xxxx}\\&\le&
7\|u_x\|_{L^\infty}\int_I\rho_{xxx}^2+
12\|u_{xx}\|_{L^\infty}\|\rho_{xx}\|_{L^2}\|\rho_{xxx}\|_{L^2}
\\&& +8\|\rho_x\|_{L^\infty}\|\rho_{xxx}\|_{L^2}\|u_{xxx}\|_{L^2}+
2\|\rho\|_{L^\infty}\|\rho_{xxx}\|_{L^2}\|u_{xxxx}\|_{L^2}. \eeqno
By Lemma \ref{non-le:4.1} and the Cauchy inequality, we get
\be\label{non-3.42} \frac{d}{dt}\int_I\rho_{xxx}^2\le
c\int_I\rho_{xxx}^2+c\int_Iu_{xxxx}^2+c.\ee
 Differentiating
(\ref{non-3.29}) with respect to $x$, we have \beq\label{non-3.43}
&&\nonumber u_{xxxx}=\rho_{xx}u_t+2\rho_xu_{xt} +\rho
u_{xxt}+(\rho_xuu_x)_x+(\rho u_x^2)_x+(\rho uu_{xx})_x\\&&\ \ \ \ \
\ \ \ \ \ \ +(P_c)_{xxx}+(\rho Q)_{xxx}. \eeq By (\ref{non-3.43}),
($A_6$) and Lemma \ref{non-le:4.1}, we have
\be\label{non-3.44} \int_Iu_{xxxx}^2 \le c\int_I\rho
u_{xxt}^2+c\int_I\rho_{xxx}^2+c. \ee By (\ref{non-3.42}),
(\ref{non-3.44}), Corollary \ref{non-cor:4.1} and the Gronwall
inequality, we get \be\label{non-3.45} \int_I\rho_{xxx}^2\le c.
 \ee
It follows from (\ref{non-3.44}), (\ref{non-3.45}) and Corollary
\ref{non-cor:4.1} that
$$
\int_{Q_T}u_{xxxx}^2\le c.
$$
A direct calculation, combining (\ref{non-1.2})$_1$,
(\ref{non-3.27}), (\ref{non-3.45}), Lemma \ref{non-le:4.1},
Corollary \ref{non-cor:4.1} and Lemma {\ref{non-le:4.2}}, implies
$$
\int_I\left(\rho_{xxt}^2+
+\rho_{tt}^2\right)+\int_{Q_T}\rho_{xtt}^2\le c.
$$
The proof of Lemma \ref{non-le:4.3} is complete. \qed

 \vspace{4mm}

The next lemma play the most important role in getting $H^4$
estimates of $u$.
\begin{Lemma}\label{non-le:4.4}Under the conditions of Theorem \ref{non-th:4.1}, it holds for any $0\le t\le T$
$$
\int_I\rho^3 u_{tt}^2+\int_{Q_T}\rho^2 u_{xtt}^2\le c.
$$
\end{Lemma}
\noindent{\it Proof.} Differentiating (\ref{non-3.22}) with respect
to $t$, we have \be\label{non-3.49} (\rho
u_{tt})_t+\rho_{tt}u_t+\rho_tu_{tt}+(\rho_tuu_x+\rho u_tu_x+\rho
uu_{xt})_t+P_{xtt}=u_{xxtt}. \ee Multiplying (\ref{non-3.49}) by
$\rho^{\gamma_2} u_{tt}$ ($\gamma_2$ is to be decided later), and
integrating the resulting equation over $I$, we have

\beqno &&\frac{1}{2}\frac{d}{dt}\int_I\rho^{\gamma_2+1}
u_{tt}^2+\int_I\rho^{\gamma_2}
u_{xtt}^2\\&=&\frac{\gamma_2-3}{2}\int_I\rho^{\gamma_2}\rho_tu_{tt}^2-\int_I[\rho_{tt}u_t
+\rho_{tt}uu_x+2\rho_tu_tu_x+2\rho_tuu_{xt}+2\rho
u_tu_{xt}+P_{xtt}](\rho^{\gamma_2}
u_{tt})\\&&-\int_I\rho^{\gamma_2+1}u_{tt}^2u_x-\int_I\rho^{\gamma_2+1}
uu_{tt}u_{xtt}-\gamma_2\int_I\rho^{\gamma_2-1}\rho_xu_{tt}u_{xtt}
\\ & \le & c\|(\sqrt{\rho})_t\|_{L^\infty}\int_I\rho^{\gamma_2+\frac{1}{2}}
u_{tt}^2+c\int_I\rho^{2\gamma_2}u_{tt}^2-\int_I\rho^{\gamma_2}
u_{tt}(\rho Q)_{xtt}-\int_I\rho^{\gamma_2} u_{tt}(P_c)_{xtt}
\\ && +c\int_I\rho^{\gamma_2+1}u_{tt}^2-\int_I\rho^{\gamma_2+1}u
u_{tt}u_{xtt}-2\gamma_2\int_I\rho^{\gamma_2-\frac{1}{2}}
u_{xtt}u_{tt}(\sqrt{\rho})_x+c
\\ & \le & c\int_I\rho^{\gamma_2+\frac{1}{2}}
u_{tt}^2+c\int_I\rho^{2\gamma_2}u_{tt}^2+\int_I\rho^{\gamma_2}
u_{xtt}(\rho
Q)_{tt}+\gamma_2\int_I\rho^{\gamma_2-1}\rho_xu_{tt}(\rho Q)_{tt} \\
&& +c\|\rho_{xtt}\|_{L^2}^2 +\frac{1}{4}\int_I\rho^{\gamma_2}
u_{xtt}^2+c\int_I\rho^{\gamma_2+2}u_{tt}^2+c\int_I\rho^{\gamma_2-1}u_{tt}^2+c\\&\le&c\int_I\rho^{\gamma_2-1}
u_{tt}^2+c\int_I\rho^{2\gamma_2}u_{tt}^2+\frac{1}{2}\int_I\rho^{\gamma_2}
u_{xtt}^2+c\int_I\rho^{\gamma_2}\left|(\rho Q)_{tt}\right|^2
\\ && +c\int_I\rho^{2\gamma_2-2}u_{tt}^2+c\int_I\rho\left|(\rho
Q)_{tt}\right|^2 +c\|\rho_{xtt}\|_{L^2}^2+c. \eeqno Here, we have
used integration by parts, the Cauchy inequality, ($A_2$), ($A_3$),
Lemma \ref{non-le:2.2}, Lemma \ref{non-le:4.1}, Lemma
\ref{non-le:4.2} and Lemma \ref{non-le:4.3}.

After the third term of the right side is absorbed by the left, we
have \beq\label{non-r3.53}
\nonumber\frac{d}{dt}\int_I\rho^{\gamma_2+1}
u_{tt}^2+\int_I\rho^{\gamma_2} u_{xtt}^2 & \le & \nonumber
c\int_I\rho^{\gamma_2-1}
u_{tt}^2+c\int_I\rho^{2\gamma_2}u_{tt}^2+c\int_I\rho^{\gamma_2}\left|(\rho
Q)_{tt}\right|^2\\&&+c\int_I\rho^{2\gamma_2-2}u_{tt}^2+c\int_I\rho\left|(\rho
Q)_{tt}\right|^2+c\|\rho_{xtt}\|_{L^2}^2+c.
 \eeq
By Lemma \ref{non-le:4.2}, we know $\int_{Q_T}\rho u_{tt}^2\le c$.
This implies that the minimum of $\gamma_2$ we should take in
(\ref{non-r3.53}) is $2$. Substituting $\gamma_2=2$ into
(\ref{non-r3.53}), we have
\be\label{non-3.50}\frac{d}{dt}\int_I\rho^3
u_{tt}^2+\int_I\rho^2u_{xtt}^2\le c\int_I\rho
u_{tt}^2+c\int_I\rho\left|(\rho
Q)_{tt}\right|^2+c\int_I\rho_{xtt}^2+c.\ee We are going to estimate
$\int_I\rho\left|(\rho Q)_{tt}\right|^2$. Using Lemma
\ref{non-le:4.1}, ($A_4$) and Lemma \ref{non-le:4.3}, we
have\beq\label{non-3.51}\nonumber \int_I\rho\left|(\rho
Q)_{tt}\right|^2&=&\int_I\rho\left|\rho_{tt}
Q+2\rho_tQ^\prime\theta_t+\rho Q^{\prime\prime}\theta_t^2+\rho
Q^\prime\theta_{tt}\right|^2\\&\le&c+c\|\theta_t\|_{L^\infty}^2+c\int_I\rho^3\theta_{tt}^2.
\eeq Substituting (\ref{non-3.51}) into (\ref{non-3.50}),
integrating the resulting inequality over $(0,t)$, and using Lemma
\ref{non-le:4.1}, Lemma \ref{non-le:4.2} and Lemma \ref{non-le:4.3},
we get \be\label{non-3.52}\int_I\rho^3
u_{tt}^2+\int_0^t\int_I\rho^2u_{xtt}^2\le\int_I\rho^3 u_{tt}^2(x,
0)+c. \ee Using (\ref{non-4.1}), (\ref{non-r3.1}), (\ref{non-r3.2}),
(\ref{non-4.6}), (\ref{non-3.7}), (\ref{non-3.22}) and
(\ref{non-r3.1})$_3$, we have
$$
\int_I\rho^3u_{tt}^2(x,0)\le c,
$$
which combining (\ref{non-3.52}) completes the proof. \qed

\begin{Lemma}\label{non-le:4.5}Under the conditions of Theorem \ref{non-th:4.1}, it holds for any $0\le t\le T$
\beqno \int_I\rho
u_{xxt}^2+\int_{Q_T}\left(u_{xxxt}^2+\rho\theta_{xxt}^2\right)\le c.
 \eeqno
\end{Lemma}
\noindent{\it Proof.} By (\ref{non-3.22}), we have \beqno \int_I\rho
u_{xxt}^2&\le&c\int_I\rho^3u_{tt}^2+c\int_I\rho\rho_t^2u_t^2+c\int_I\rho\rho_t^2u^2u_x^2+\int_I\rho^3
u_t^2u_x^2 \\ && +c\int_I\rho^3u^2u_{xt}^2+c\int_I\rho|(\rho
Q)_{xt}|^2+c\int_I\rho|(P_c)_{xt}|^2\\&\le&c, \eeqno where we have
used ($A_2$)-($A_4$), Lemma \ref{non-le:2.2}, Lemma
\ref{non-le:4.1}, Lemma \ref{non-le:4.2} and Lemma \ref{non-le:4.4}.

It follows from (\ref{non-3.36}) and ($A_5$) that \beqno
\int_{Q_T}\rho\theta_{xxt}^2&\le&c\int_{Q_T}\rho|\kappa^\prime|^2\theta_t^2\theta_{xx}^2+c\int_{Q_T}\rho|\kappa^{\prime\prime}|^2\theta_t^2\theta_x^4
+c\int_{Q_T}\rho|\kappa^\prime|^2\theta_x^2\theta_{xt}^2 \\ & &
+c\int_{Q_T}\rho^3
|Q^\prime|^2\theta_{tt}^2+c\int_{Q_T}\rho^3|Q^{\prime\prime}|^2\theta_t^4+c\int_{Q_T}\rho\rho_t^2|Q^\prime|^2\theta_t^2
\\ && +c\int_{Q_T}\rho\left|(\rho
uQ^\prime\theta_x)_t\right|^2+c\int_{Q_T}\rho\left|(\rho\theta
Q^\prime u_x)_t\right|^2+c\int_{Q_T}\rho u_x^2u_{xt}^2,\eeqno which,
combining ($A_4$), ($A_5$) and Lemma \ref{non-le:4.1}, gives
\beq\label{non-3.55}
\nonumber\int_{Q_T}\rho\theta_{xxt}^2&\le&c\int_{Q_T}\theta_{xt}^2+c\int_{Q_T}\rho^3\theta_{tt}^2+c\int_0^T\|\theta_t\|_{L^\infty}^2
+\int_{Q_T}u_{xt}^2+c\\&\le&c. \eeq Differentiating (\ref{non-3.29})
with respect to $t$, we get \beqno
u_{xxxt}&=&2(\sqrt{\rho})_x\sqrt{\rho}u_{tt}+\rho
u_{xtt}+\rho_{xt}u_t+\rho_tu_{xt}
 +\rho_{xt}uu_x+\rho_tu_x^2+\rho_tuu_{xx}+\rho_xu_tu_x\\&&+2\rho
u_{xt}u_x+\rho u_tu_{xx}+\rho_xuu_{xt}+\rho uu_{xxt}+(\rho
Q)_{xxt}+(P_c)_{xxt}. \eeqno This, together with (\ref{non-3.55}),
($A_6$), Lemma \ref{non-le:2.2}, Lemma \ref{non-le:4.1}, Lemma
\ref{non-le:4.2}, Lemma \ref{non-le:4.3}, Lemma \ref{non-le:4.4} and
Corollary \ref{non-cor:4.1}, implies \beqno \int_{Q_T}u_{xxxt}^2 \le
c+c\int_{Q_T}\left(\rho\theta_{xxt}^2+\theta_{xt}^2\right)+c\int_0^T\|\theta_t\|_{L^\infty}^2+c\int_{Q_T}\rho_{xxt}^2
\le c.
\eeqno This completes the proof.\qed\\[1.5mm]

From (\ref{non-3.44}), Lemma \ref{non-le:4.3} and Lemma
\ref{non-le:4.5}, we get the following corollary immediately.
\begin{cor}\label{non-cor:4.2}Under the conditions of Theorem \ref{non-th:4.1}, it holds for any $0\le t\le T$
$$\int_Iu_{xxxx}^2\le c.$$
\end{cor}

\begin{cor}\label{non-cor:4.3}Under the conditions of Theorem \ref{non-th:4.1}, it holds
$$\int_{Q_T}\theta_{xxxx}^2\le c.$$
\end{cor}
\noindent{\it Proof.} Differentiating (\ref{non-3.39}) with respect
to $x$, we have \beqno
&&\kappa\theta_{xxxx}=-4\kappa^\prime\theta_x\theta_{xxx}-3\kappa^\prime\theta_{xx}^2-3\kappa^{\prime\prime}\theta_x^2\theta_{xx}
-\left(\kappa^{\prime\prime}\theta_x^3\right)_x-2\left(u_xu_{xx}\right)_x+\left(\rho
Q^\prime\theta_{xt}\right)_x\\&&\ \ \ \ \ \ \ \ \ \ \ \
+\left(\rho_xQ^\prime\theta_t\right)_x+\left(\rho
Q^{\prime\prime}\theta_x\theta_t\right)_x +(\rho u
Q^\prime\theta_x)_{xx}+(\rho\theta Q^\prime u_x)_{xx}. \eeqno This,
combining ($A_5$), ($A_6$), Lemma \ref{non-le:3.3}, Lemma
\ref{non-le:4.1} and Lemma \ref{non-le:4.5}, implies \beqno
\int_{Q_T}\theta_{xxxx}^2 \le
c+c\int_{Q_T}\left(\rho\theta_{xxt}^2+\theta_{xt}^2\right)+c\int_0^T\|\theta_t\|_{L^\infty}^2
+c\int_{Q_T}|\kappa^{\prime\prime\prime}|^2\theta_x^8 \le c. \eeqno
This proves Corollary 4.3. \qed

\begin{Lemma}\label{non-le:4.7}Under the conditions of Theorem \ref{non-th:4.1}, it holds for any $0\le t\le T$
$$
\int_I\rho_{xxxx}^2+\int_{Q_T}u_{xxxxx}^2\le c.
$$

\end{Lemma}
\noindent{\it Proof.} Differentiating (\ref{non-3.41}) with respect
to $x$, multiplying the resulting equation by $2\rho_{xxxx}$,
integrating over $I$, and using integration by parts, Lemma
\ref{non-le:4.1}, Lemma \ref{non-le:4.3}, Corollary
\ref{non-cor:4.2} and the Cauchy inequality, we get
\beq\label{non-3.56}\nonumber
\frac{d}{dt}\int_I\rho_{xxxx}^2&=&-9\int_I\rho_{xxxx}^2u_x
-20\int_I\rho_{xxx}\rho_{xxxx}u_{xx}-20\int_I\rho_{xx}\rho_{xxxx}u_{xxx}
\\&&\nonumber -10\int_I\rho_x\rho_{xxxx}u_{xxxx}-2\int_I\rho\rho_{xxxx}u_{xxxxx}
\\&\le&c\int_I\rho_{xxxx}^2+c\int_Iu_{xxxxx}^2+c.\eeq
 Now we estimate the second term of the right-hand side of
(\ref{non-3.56}).

Differentiating (\ref{non-3.43}) with respect to $x$, we have \beqno
u_{xxxxx}&=&\rho_{xxx}u_t+3\rho_{xx}u_{xt}+3\rho_xu_{xxt}+\rho
u_{xxxt}+(\rho_xuu_x)_{xx}+(\rho u_x^2)_{xx}
\\&& +(\rho
uu_{xx})_{xx}+(\rho Q)_{xxxx}+(P_c)_{xxxx}. \eeqno This, combining $(A_6)$,
Lemma \ref{non-le:2.2}, Lemma \ref{non-le:4.1}, Lemma
\ref{non-le:4.3} and Corollary \ref{non-cor:4.2}, concludes
\beq\label{non-3.57}\int_I u_{xxxxx}^2\le
c\int_I u_{xxt}^2+c\int_Iu_{xxxt}^2
+c\int_I\rho_{xxxx}^2+c\int_I\theta_{xxxx}^2+c. \eeq Substituting
(\ref{non-3.57}) into (\ref{non-3.56}), and using Corollary
\ref{non-cor:4.1}, Corollary \ref{non-cor:4.3}, Lemma
\ref{non-le:4.5} and the Gronwall inequality, we obtain
\be\label{non-3.58} \int_I\rho_{xxxx}^2\le c. \ee It follows from
(\ref{non-3.57}), (\ref{non-3.58}), Corollary \ref{non-cor:4.1},
Corollary \ref{non-cor:4.3} and Lemma \ref{non-le:4.5} that
$$\int_{Q_T}u_{xxxxx}^2\le c.
$$
This completes the proof of Lemma 4.6. \qed

\begin{cor}\label{non-cor:4.4}Under the conditions of Theorem \ref{non-th:4.1}, it holds for any $0\le t\le T$
$$
\int_I\left(\rho_{xtt}^2+\rho_{xxxt}^2
\right)+\int_{Q_T}\left(\rho_{ttt}^2+\rho_{xxtt}^2 \right)\le c.
$$
\end{cor}
Here we have used the following inequality when we get the upper
bound of $\rho_{ttt}$:
$$
\rho_x^2u_{tt}^2=2\left[(\sqrt{\rho})_x\sqrt{\rho}\right]^2u_{tt}^2\leq
c\rho u_{tt}^2.
$$
From the above estimates, we get \beq\label{non-3.60}\nonumber
&&\|(\sqrt{\rho})_x\|_{L^\infty}+\|(\sqrt{\rho})_t\|_{L^\infty}+\|\rho\|_{H^4}
+\|\rho_t\|_{H^3}
 \displaystyle +\|\rho_{tt}\|_{H^1}+\|u\|_{H^4}
\\ && \nonumber +\|u_t\|_{H^1}+\|\rho^\frac{3}{2}u_{tt}\|_{L^2}+\|\sqrt{\rho}u_{xxt}\|_{L^2}+\|\theta\|_{H^3}+\|\sqrt{\rho}\theta_t\|_{L^2}+\|\rho\theta_{xt}\|_{L^2}
\\ && \nonumber +\int_{Q_T}\left(\rho^2 u_{xtt}^2+\rho
u_{tt}^2+u_{xxt}^2+u_{xxxt}^2+u_{xxxxx}^2\right)\\ &&
+\int_{Q_T}\left(\rho_{ttt}^2
+\rho_{xxtt}^2+\theta_{xxxx}^2+\theta_t^2+\theta_{xt}^2+\rho\theta_{xxt}^2+\rho^3\theta_{tt}^2\right)\le
c. \eeq From (\ref{non-3.60}) and (\ref{non-r3.55}), we get \beqno
&&\|\rho\|_{H^4}+\|\rho_t\|_{H^3}+\|\rho_{tt}\|_{H^1}+\|u\|_{H^4}
+\|u_t\|_{H^2}+\|u_{tt}\|_{L^2}+\|\theta\|_{H^3}+\|\theta_t\|_{H^1}\\&&
+\int_{Q_T}\left(u_{xtt}^2+u_{xxxt}^2+u_{xxxxx}^2\right)
+\int_{Q_T}\left(\rho_{ttt}^2
+\rho_{xxtt}^2+\theta_{xxxx}^2+\theta_{xxt}^2+\theta_{tt}^2\right)
\le c(\delta), \eeqno where $c(\delta)$ is a positive constant, and
may depend on $\delta$.

The proof of Theorem \ref{non-th:4.1} is complete. \qed\\[1.5mm]

{\noindent\bf Proof of Theorem \ref{non-th:1.1}:}

Consider (\ref{non-1.2})-(\ref{non-1.4}) with initial data replaced
by ($\rho_0^\delta$, $u_0^{\delta}$, $\theta_0^\delta$), we obtain
from Theorem \ref{non-th:4.1} that there exists a unique solution
($\rho^{\delta}$, $u^{\delta}$, $\theta^{\delta}$) such that
(\ref{non-3.60}) and (\ref{non-r3.55}) are valid when we replace
($\rho,$ $u$, $\theta$) by ($\rho^{\delta}$, $u^{\delta}$,
$\theta^{\delta}$). With this estimates uniform for $\delta$, we
take $\delta\rightarrow0^+$ ( take subsequence if necessary) to get a
solution to (\ref{non-1.2})-(\ref{non-1.4}) still denoted by
($\rho$, $u$, $\theta$). By the lower semi-continuity of the norms,
we have
 \beqno
&&\|(\sqrt{\rho})_x\|_{L^\infty}+\|(\sqrt{\rho})_t\|_{L^\infty}+\|\rho\|_{H^4}
+\|\rho_t\|_{H^3}
 \displaystyle +\|\rho_{tt}\|_{H^1}+\|u\|_{H^4}
+\|u_t\|_{H^1}
\\ && +\|\sqrt{\rho}u_{xxt}\|_{L^2}+\|\theta\|_{H^3}+\|\sqrt{\rho}\theta_t\|_{L^2}+\|\rho\theta_{xt}\|_{L^2}
+\int_{Q_T}\left({u}_{xxt}^2+{u}_{xxxt}^2+{u}_{xxxxx}^2\right)\\
&& +\int_{Q_T}\left({\rho}_{ttt}^2
+{\rho}_{xxtt}^2+{\theta}_{xxxx}^2+{\theta}_t^2+{\theta}_{xt}^2\right)\le
c, \eeqno  which proves the existence of the solutions as in Theorem
\ref{non-th:1.1}. The uniqueness of the solutions can be proved by
the standard method like in \cite{cho-Kim: perfect gas}, we omit it for brevity. The
proof of Theorem \ref{non-th:1.1} is complete.\qed

\vspace{6mm}

 \noindent{\bf\small Acknowledgment.} {\small

The first author was supported by the National Basic Research
Program of China (973 Program) No. 2010CB808002, and by the National
Natural Science Foundation of China No. 11071086. The second author
was supported by the National Natural Science Foundation of China
$\#$10625105, $\#$11071093, the PhD specialized grant of the
Ministry of Education of China $\#$20100144110001 and the
self-determined research funds of CCNU from the colleges'basic
research and operation of MOE. }

\end{document}